\documentclass[12pt,a4paper]{amsart}

\usepackage{{amsmath,latexsym,amssymb,amscd}}
\usepackage{fullpage}                               
\title{Families of Monads and Instantons from \\[6pt] a Noncommutative ADHM Construction}

\author{Simon Brain}
\address{\flushleft SISSA, Via Beirut 2-4, 34151 Trieste, Italy} \email{brain@sissa.it}
\author{Giovanni Landi}
\address{\flushleft Dipartimento di Matematica e
Informatica, Universit\`{a} di Trieste, Via A. Valerio 12/1, 34127
Trieste, Italy, and INFN, Sezione di Trieste, Trieste, Italy}
\email{landi@univ.trieste.it}
\thanks{Partially supported by:}
\thanks{the `Italian project Cofin06 - Noncommutative geometry, quantum groups and applications'}
\date{v1: 7 January 2009; v2: 6 August 2010}

\newtheorem{thm}{Theorem}[section]

\newtheorem{lem}[thm]{Lemma}
\newtheorem{prop}[thm]{Proposition}
\newtheorem{rem}[thm]{Remark}

\theoremstyle{definition}
\newtheorem{defn}[thm]{Definition}

\newcommand{\ii}{\mathrm{i}}

\newcommand{\pp}{{\sf p}}
\newcommand{\qp}{{\sf q}}
\newcommand{\Pp}{{\sf P}}
\newcommand{\Qp}{{\sf Q}}

\newcommand{\E}{\mathcal{E}}

\newcommand{\A}{\mathcal{A}}

\newcommand{\la}{\langle}
\newcommand{\ra}{\rangle}

\newcommand{\tr}{\triangleright}
\newcommand{\tl}{\triangleleft}

\newcommand{\n}{\nabla}
\newcommand{\ep}{\epsilon}

\newcommand{\M}{\textup{M}}
\newcommand{\tmM}{\widetilde{\mathcal{M}}}
\newcommand{\SL}{\textup{SL}}
\newcommand{\Sp}{\textup{Sp}}
\newcommand{\GL}{\textup{GL}}
\newcommand{\U}{\textup{U}}
\newcommand{\V}{\textup{V}}

\newcommand{\SU}{\textup{SU}}

\newcommand{\CP}{\mathbb{C} \mathbb{P}}

\newcommand{\HH}{\mathbb{H}}
\newcommand{\ZZ}{\mathbb{Z}}
\newcommand{\C}{\mathbb{C}}

\newcommand{\RR}{\mathbb{R}}

\newcommand{\End}{\textup{End}}

\newcommand{\id}{\textup{id}}
\newcommand{\D}{\textup{d}}

\begin{document}

%%\title[Noncommutative Families of Monads and Instantons]{Families of Monads and Instantons \\ from a Noncommutative ADHM Construction}
%%%\date{31 January 2009}
%%%\date{16 July 2009}
%%%

%%
%%\author{Simon Brain}
%%%\address{\flushleft SISSA, Via Beirut 2-4, 34151 Trieste, Italy} \email{brain@sissa.it}
%%\address{SISSA, Via Beirut 2-4, 34151 Trieste, Italy} \email{brain@sissa.it}
%%\author{Giovanni Landi}
%%%\address{\flushleft Dipartimento di Matematica e
%%%Informatica, Universit\`{a} di Trieste, Via A. Valerio 12/1, 34127
%%%Trieste, Italy, and INFN, Sezione di Trieste, Trieste, Italy}
%%\address{Dipartimento di Matematica e
%%Informatica, Universit\`{a} di Trieste, Via A. Valerio 12/1, 34127
%%Trieste, Italy, and INFN, Sezione di Trieste, Trieste, Italy}

%%\email{landi@univ.trieste.it}
%%\copyrightinfo{Simon Brain and Giovanni Landi}

%%\thanks{Partially supported by:}
%%\thanks{the `Italian project Cofin06 - Noncommutative geometry, quantum groups and applications'}

\dedicatory{Dedicated to Alain Connes}

\subjclass[2000]{Primary 58B34; Secondary 14D21, 81T13}

\maketitle

\begin{abstract}
We give a $\theta$-deformed version of the ADHM construction of
instantons with arbitrary topological charge on the sphere $S^4$.
Classically, the instanton gauge fields are constructed from
suitable monad data; we show that in the deformed case the set of
monads is itself a noncommutative space. We use these monads to
construct noncommutative `families' of instantons (i.e.
noncommutative families of anti-self-dual connections) on the
deformed sphere $S^4_\theta$. We also compute the topological charge
of each of the families. Finally we discuss what
it means for such families to be gauge equivalent. \\
\end{abstract}

\tableofcontents

%\newpage

\section{Introduction}

The purpose of the present article is to generalise the ADHM method for constructing instantons
on the four-sphere $S^4$ to the framework of noncommutative geometry, by giving a
construction of instantons on the noncommutative four-sphere $S^4_\theta$ of \cite{cl:id}.

Instantons arise in physics as anti-self-dual solutions of the
Yang-Mills equations. Mathematically they are connections with
anti-self-dual curvature on smooth $G$-bundles over a
four-dimensional compact manifold. Since the very beginning they
have been of central importance for both disciplines, an importance
that has only grown over the years.

Of particular interest are instantons on $\SU(2)$-bundles over the Euclidean four-sphere $S^4$.
Thanks to the ADHM method of \cite{adhm:ci}, the full solution to the
problem of constructing such instantons on $S^4$ has long been known and,
as a consequence, the moduli space $\mathcal{M}_k$ of
%$\SU(2)$
instantons with topological charge equal to $k$ is known to be a
manifold of dimension $8k-3$. Starting with a trivial vector bundle
over $S^4$, the ADHM strategy is to construct an orthogonal
projection to some (non-trivial) sub-bundle $E$ in such a way that
the projection of the trivial connection to $E$ has anti-self-dual
curvature.

The geometric ingredient which implements the classical ADHM
construction is the Penrose twistor fibration $\CP^3\to S^4$. The
total space $\CP^3$ of the fibration is called the twistor space of
$S^4$ and may be thought of as the bundle of projective spinors over
$S^4$ (although it has its origins elsewhere \cite{pr:sst}). The
pull-back of an instanton bundle along this fibration is a
holomorphic vector bundle over $\CP^3$ equipped with a set of
reality conditions which identify it as such a pull-back
\cite{rw:sdgf}. In this way, the construction of instantons is
equivalent to the construction of holomorphic bundles over twistor
space.

Using powerful results from algebraic geometry, one
gives an explicit description of all relevant holomorphic vector
bundles over a complex projective space (\cite{gh:vbps,wb:mon}, {\em
cf}. also \cite{oss:vb}). Each of them arises as the cohomology of a
{\em monad}: a suitable complex of vector bundles
$$0\to \mathcal{A} \xrightarrow{\sigma} \mathcal{B}
\xrightarrow{\tau} \mathcal{C} \to 0$$ such that $\sigma$ is
injective and $\tau$ is surjective. The ADHM construction tells us
how to convert a given monad into an orthogonal projection of vector
bundles as described above and guarantees that the resulting
connection has anti-self-dual curvature.

Following the general strategy of the classical case, our goal
%in what follows
is to give a deformed version of the ADHM method and hence a
construction of %$\SU(2)$
instantons on the noncommutative four-sphere $S^4_\theta$. The
techniques involved lend themselves rather neatly to the framework
of noncommutative geometry; the construction of vector bundles and
connections by orthogonal projection is particularly natural in
light of the Serre-Swan theorem \cite{dvg:pvf}, which trades vector
bundles for finitely generated projective modules.

The paper is organised as follows. Sect.~\ref{section twistor
fibration} reviews the noncommutative spaces in question, namely the
$\theta$-deformed versions of the four-sphere $S^4_\theta$ and its
twistor space $\CP^3_\theta$. We recall also the construction of the
basic instanton and the principal bundle on which it is defined, as
well as the details of the noncommutative twistor fibration.
Sect.~\ref{section noncommutative conformal transformations} recalls
the construction of the quantum group $\SL_\theta(2,\HH)$ of
conformal transformations of $S^4_\theta$ and the quantum subgroup
$\Sp_\theta(2)$ of isometries. The main purpose of these two
sections is to gather together into one place the relevant
contributions from \cite{cl:id,lvs:pfns,lvs:nitcs,lprs:ncfi,bm:qtt}
and to establish notation; in doing so we also make some novel
improvements to previous versions. Sect.~\ref{section adhm
construction} presents the deformed ADHM construction itself. We
show that in the deformed case the set of all monads is
parameterised by a collection of noncommutative spaces
$\tmM_{\theta;k}$ indexed by $k$ a positive integer. We use each of
these spaces to construct a noncommutative `family' of instantons
whose topological charge we show to be equal to $k$. Finally in
Sect.~\ref{section gauge theory} we discuss what it means for
families of instantons to be gauge equivalent. In particular, we
show that the quantum symmetries of the sphere $S^4_\theta$ generate
gauge degrees of freedom, a feature which is a consequence of the
noncommutativity and is not present in the classical construction. 
For further discussion in this direction we refer to \cite{bl:adhm}.

\section{The Twistor Fibration}\label{section twistor fibration}

The use of the twistor fibration in the ADHM construction is crucial: this fibration
captures in its geometry the very nature of the
anti-self-duality equations,
%{\color{red} \cite{pr:sst,mw:isdtt}}
with the result that an instanton bundle is reinterpreted {\em via}
pull-back in terms of holomorphic data on twistor space
\cite{rw:sdgf} ({\em cf}. also \cite{aw:iag}). In particular, this
means that twistor space plays the role of an `auxiliary space' on
which the ADHM construction takes place, before passing back down to
the base space $S^4$ (we refer to \cite{mw:isdtt} for more on the
ADHM construction from a twistor perspective).

We start by recalling the details of the algebra
inclusion $\A(S^4_\theta) \hookrightarrow \A(S^7_\theta)$ as a
noncommutative principal bundle with undeformed structure group
$\SU(2)$; associated to this principal bundle there is in particular a basic
instanton bundle \cite{lvs:pfns}. Next we give a description of the
noncommutative twistor space in terms of its coordinate algebra
$\A(\CP^3_\theta)$, as well as a dualised description of the twistor
fibration, now appearing \cite{bm:qtt} as an algebra inclusion
$\A(S^4_\theta)\hookrightarrow \A(\CP^3_\theta)$.

\subsection{The noncommutative Hopf fibration} \label{section noncommutative hopf fibration}
With $\lambda=\exp{(2\pi\ii \theta)}$ the deformation parameter, the
coordinate algebra $\A(S^4_\theta)$ of the noncommutative
four-sphere $S^4_\theta$ is the $*$-algebra generated by a central
real element $x$ and elements $\alpha$, $\beta$, $\alpha^*$,
$\beta^*$, modulo the relations \begin{equation} \label{eqn four
sphere relations} \alpha \beta= \lambda \beta \alpha, \quad \alpha^*
\beta^* = \lambda \beta^* \alpha ^*, \quad \beta^* \alpha = \lambda
\alpha \beta^*, \quad \beta \alpha^* = \lambda \alpha^*
\beta,\end{equation} together with the sphere relation
\begin{equation}\label{eqn foursphere relation}\alpha^*\alpha + \beta^* \beta + x^2=1.
\end{equation}
Similarly, the coordinate algebra of the noncommutative seven-sphere
$\A(S^7_\theta)$ is generated as a $*$-algebra by the elements $\{
z_j, z_j^* ~|~ j=1,\ldots,4 \}$ and is subject to the commutation
relations \begin{equation} \label{eqn nc params} z_j
z_l=\eta_{jl}z_l z_j, \quad z_jz_l^*=\eta_{lj} z_l^*z_j, \quad
z^*_jz^*_l=\eta_{jl}z^*_lz^*_j,\end{equation} as well as the sphere
relation \begin{equation}\label{eqn sevensphere
relation}z^*_1z_1+z^*_2z_2+z^*_3z_3+z^*_4z_4=1.\end{equation}
Compatibility with the $\SU(2)$ principal bundle structure requires
the deformation matrix $(\eta_{jk})$ be given by
\begin{equation}\label{eqn eta matrix}(\eta_{jk})=\begin{pmatrix} 1 & 1
& \bar\mu & \mu \\ 1 & 1 & \mu & \bar\mu \\ \mu & \bar \mu & 1 & 1 \\
\bar\mu & \mu & 1 & 1 \end{pmatrix}, \qquad \mu=\exp{(\ii \pi
\theta)}.\end{equation}
The values of the deformation parameters
$\lambda$, $\mu$ are precisely those which allow an embedding of the
classical group $\SU(2)$ into the group
$\textup{Aut}\,\A(S^7_\theta)$. We denote by $\A(\C^4_\theta)$ the algebra
generated by the $\{z_j$, $z^*_j\}$ subject to the relations (\ref{eqn
nc params}); the quotient by the additional sphere relation yields
the algebra $\A(S^7_\theta)$. The algebra inclusion $\A(S^4_\theta)
\hookrightarrow \A(S^7_\theta)$ is given explicitly by
\begin{eqnarray} \label{eqn algebra inclusion}
&\alpha=2(z_1z^*_3 + z_2^*z_4), \quad
\beta=2(z_2z_3^* - z_1^*z_4), \quad x= z_1z^*_1 + z_2z^*_2 -
z_3z_3^* - z_4z_4^*.\end{eqnarray} One easily verifies that for
the right $\SU(2)$-action on $\A(S^7_\theta)$ given on generators by
\begin{equation} \label{eqn su(2) action}
(z_1,z_2^*,z_3,z_4^*)\mapsto(z_1,z_2^*,z_3,z_4^*)\begin{pmatrix}w&0\\0&w\end{pmatrix},\qquad
w=\begin{pmatrix}w^1 & -\bar w^2\\ w^2& \bar w^1\end{pmatrix}\in
\SU(2),\end{equation} the invariant subalgebra is generated as
expected by $\alpha$, $\beta$, $x$ and their conjugates, so one
indeed has
$$\textup{Inv}_{\SU(2)}\A(S^7_\theta)=\A(S^4_\theta).$$ When $\theta=0$ we recover the
usual algebras of functions on the classical spheres $S^4$ and
$S^7$. The inclusion $\A(S^4) \hookrightarrow \A(S^7)$ is just a
dualised description of the standard $\SU(2)$ Hopf fibration $S^7
\to S^4$.

These noncommutative spheres have canonical differential calculi
arising as deformations of the classical ones. Explicitly, one has a
first order differential calculus $\Omega^1(S^7_\theta)$ on
$\A(S^7_\theta)$ spanned as an $\A(S^7_\theta)$-bimodule by $\{ \D
z_j$, $\D z_j^*$, $j=1,\ldots,4\}$, subject to the relations $$z_i
\D z_j = \eta_{ij} \D z_j z_i, \qquad z_i \D z_j^* =\eta_{ji} \D
z_j^* z_i,$$ with $\eta_{ij}$ as before. One also has relations
$$\D z_i \D z_j + \eta_{ij} \D z_j \D z_i=0, \qquad \D z_i \D
z_j^* + \eta_{ji} \D z_j^* \D z_i=0,$$ allowing one to extend the
first order calculus to a differential graded algebra
$\Omega(S^7_\theta)$. There is a unique differential $\D$ on
$\Omega(S^7_\theta)$ such that $\D: z_j \mapsto \D z_j$.
Furthermore, $\Omega(S^7_\theta)$ has an involution given by the
graded extension of the map $z_j \mapsto z_j^*$. The story is
similar for the four-sphere, in that the differential graded algebra
$\Omega(S^4_\theta)$ is generated in degree
%zero by $\alpha$, $\alpha^*$, $\beta$, $\beta^*$, $x$ and in degree
one by $\D \alpha$, $\D \alpha^*$, $\D \beta$, $\D \beta^*$, $\D x$,
subject to the relations
$$\alpha \D \beta=\lambda (\D
\beta)\alpha, \qquad \beta^* \D \alpha = \lambda (\D \alpha)
\beta^*,$$ $$\D \alpha \D \beta + \lambda \D \beta \D \alpha=0,
\qquad \D \beta^* \D \alpha + \lambda \D \alpha \D \beta^*=0 .
$$ The above
are the same as the relations (\ref{eqn four sphere relations}) or  \eqref{eqn nc params} but
with $\D$ inserted. As vector spaces, the graded components
$\Omega^k(S^7_\theta)$ and $\Omega^k(S^4_\theta)$ of $k$-forms on
the noncommutative spheres are identical to their classical
counterparts, although the algebra relations between forms are
twisted. In particular this means that the Hodge $*$-operator on $S^4_\theta$,
$$*_\theta:\Omega^k(S^4_\theta) \to \Omega^{4-k}(S^4_\theta) , $$ is
defined by the same formula as it is classically. One still has that
$*_\theta^2=1$, whence there is a direct sum decomposition of
two-forms
$$\Omega^2(S^4_\theta) = \Omega^2_+(S^4_\theta) \oplus \Omega^2_-(S^4_\theta),$$
with $\Omega^2_{\pm}(S^4_\theta):=\{\omega \in \Omega^2(S^4_\theta)
~|~ *_\theta \omega=\pm \omega\}$ the spaces of self-dual and
anti-self-dual two-forms.

\subsection{The basic instanton} \label{section basic instanton}
Amongst the nice properties of the classical
Hopf fibration is that its canonical connection is an
anti-instanton: its curvature is a self-dual two-form with values in
the Lie algebra $\mathfrak{su}(2)$ of the structure group. This
property holds also in the noncommutative case, giving a simple
example of a noncommutative instanton. It has an elegant description \cite{lvs:pfns}
in terms of the function algebras $\A(S^7_\theta)$, $\A(S^4_\theta)$
as follows. One takes the pair of elements of the right
$\A(S^7_\theta)$-module $\A(S^7_\theta)^4:=\C^4 \otimes
\A(S^7_\theta)$ given by
$$|\psi_1\ra=\begin{pmatrix}z_1 & z_2 & z_3 & z_4\end{pmatrix}^{\textup{t}}, \qquad
|\psi_2\ra=\begin{pmatrix}-z_2^* & z_1^* &
-z_4^* & z^*_3\end{pmatrix}^{\textup{t}}.$$
With the natural
Hermitian structure on $\A(S^7_\theta)^4$ given by $\la \xi |
\eta\ra = \sum_i \xi^*_i \eta_i$, one sees that
$\la\psi_j | \psi_l \ra = \delta_{jl}$. It is convenient to
introduce the matrix-valued function $\Psi$ on $S^7_\theta$
given by \begin{equation}\label{eqn partial isometry}\Psi=\begin{pmatrix}
|\psi_1\ra & |\psi_2 \ra\end{pmatrix}=\begin{pmatrix} z_1 & z_2 & z_3 & z_4 \\
-z_2^* & z_1^* & -z_4^* & z^*_3
\end{pmatrix}^{\textup{t}}.\end{equation} From
orthonormality of the columns one has that $\Psi^*\Psi=1$ and hence
the matrix
\begin{equation} \label{eqn basic instanton projector}
\qp :=\Psi\Psi^*=\frac{1}{2}\begin{pmatrix} 1+x & 0 & \alpha & -\bar
\mu\, \beta^* \\ 0 & 1+x & \beta & \mu\, \alpha^* \\ \alpha^* &
\beta^* & 1-x & 0 \\ -\mu\, \beta & \bar\mu\, \alpha & 0 & 1-x
\end{pmatrix}\end{equation} is a self-adjoint idempotent of rank two, {\em i.e.}
$\qp^*=\qp=\qp^2$ and $\textup{Tr}\, \qp=2$. The action (\ref{eqn su(2)
action}) of $\SU(2)$ on $\A(S^7_\theta)$ now takes the form
$$\Psi \mapsto \Psi w, \qquad w \in \SU(2),$$ from which the
$\SU(2)$-invariance of the entries of $\qp$ is immediately deduced. We
may also write the commutation relations of $\A(S^7_\theta)$ in the
useful form \begin{equation}\Psi_{ia} \Psi_{jb} = \eta_{ij}
\Psi_{jb} \Psi_{ia}, \qquad a,b=1,2 \quad i,j=1,2,3,4.\end{equation}

If $\rho$ is the defining representation of $\SU(2)$ on $\C^2$, the
finitely generated projective right $\A(S^4_\theta)$-module
$\E:=\qp \A(S^4_\theta)^4$ is isomorphic to the module of equivariant
maps from $\A(S^7_\theta)$ to $\C^2$,
 $$\E \cong \{\phi \in
\A(S^7_\theta) \otimes \C^2 ~|~ (w\otimes
\textup{id})\phi=(\textup{id}\otimes\rho(w^{-1}))\phi ~\text{for
all}~ w \in \SU(2)\}.$$
The module $\E$ has the role of the module
of sections of the `associated  vector bundle' $E=S^7_\theta
\times_{\SU(2)}\C^2$. With the projection
$\qp=\Psi\Psi^*$ there comes the canonical Grassmann connection
defined on the module $\E$ by
$$\n:=\qp \circ \D: \E \to \E\otimes_{\A(S^4_\theta)}
\Omega^1(S^4_\theta).$$ The curvature of $\n$ is $\n^2=\qp(\D \qp)^2$,
which may be shown to be self-dual with respect to the Hodge
operator, $$*_\theta(\qp(\D \qp)^2)=\qp(\D \qp)^2.$$ The complementary
projector $\pp=1-\qp$ yields a connection whose curvature is
anti-self-dual, $*_\theta(\pp(\D \pp)^2)=-\pp(\D \pp)^2$, and hence an instanton on the noncommutative
four-sphere, which we call the {\em basic instanton}. Noncommutative
index theory computes its `topological charge'  to be equal to $-1$.

Using the standard basis $(e_1,e_2)$ of $\C^2$,
equivariant maps are written as $\phi=\sum_a \phi_a \otimes e_a$. On them, one
has explicitly that
$$\n(\phi_a)=\D \phi_a + \sum\nolimits_b \omega_{ab}\phi_b,$$ where the
connection one-form $\omega=\omega_{ab}$ is found to be
\begin{equation} \label{eqn basic instanton one form}
\omega_{ab}=\tfrac{1}{2}
\sum\nolimits_j \left( (\Psi^*)_{aj}\D \Psi_{jb}-\D(\Psi^*)_{aj}\Psi_{jb} \right).
\end{equation}
From this
it is easy to see that $\omega_{ab}=-(\omega^*)_{ba}$ and
$\sum_a \omega_{aa}=0$, so that $\omega$ is an element of
$\Omega^1(S^7_\theta)\otimes \mathfrak{su}(2)$.

\subsection{Noncommutative twistor space} \label{section noncommutative twistor space}

It is well-known that, as a real six-dimensional manifold, the space
$\CP^3$ may be identified with the set of all $4 \times 4$ Hermitian
projector matrices of rank one: this is because each such matrix
uniquely determines and is uniquely determined by a one-dimensional
subspace of $\C^4$. Thus the coordinate algebra
$\A(\CP^3)$ of $\CP^3$ has a defining matrix of generators
\begin{equation}\label{eqn twistor matrix}
Q=\begin{pmatrix}t_1 & x_1 & x_2 & x_3\\ x_1^* & t_2 & y_3 & y_2 \\
x_2^* & y_3^* & t_3 & y_1 \\ x_3^* & y_2^* & y_1^* & t_4\end{pmatrix} ,
\end{equation}
with $t_j^*=t_j$, $j=1,\ldots,4$ and $\textup{Tr}\,Q=\sum_j t_j
=1$, as well as the relations coming from the condition $Q^2=Q$,
that is to say $\sum_j Q_{kj} Q_{jl} = Q_{kl}$. The noncommutative
twistor algebra $\A(\CP^3_\theta)$ is obtained by deforming these
relations: with deformation parameter $\lambda=\exp{(2\pi \ii
\theta)}$, one has that $t_1,\ldots, t_4$ are central, that
$$x_1x_3=\bar\lambda x_3x_1, \quad x_2x_1=\bar\lambda x_1x_2, \quad
x_2x_3=\bar\lambda x_3x_2$$ as well as the auxiliary relations
$$y_1y_2=\bar\lambda y_2y_1, \quad y_1y_3=\bar \lambda y_3y_1, \quad
y_2y_3=\bar\lambda y_3y_2, \quad
x_1(y_1,y_2,y_3)=(\bar\lambda^2y_1,\bar\lambda y_2,\lambda y_3)x_1,  $$
$$
x_2(y_1,y_2,y_3)=(\bar\lambda y_1,y_2,\lambda y_3)x_2, \quad
    x_3(y_1,y_2,y_3)=(\bar\lambda y_1,\bar\lambda y_2,y_3)x_3,
$$
and similar relations obtained by taking the adjoint under $*$ of
those above (we refer to \cite{bm:qtt} for further details). To
proceed further it is useful to note that classically $\CP^3$ is the
quotient of the sphere $S^7$ by the action of the diagonal $\U(1)$
subgroup of $\SU(2)$. This remains true in the noncommutative case
and one identifies the generators of $\A(\CP^3_\theta)$ as
\begin{equation} \label{eqn twistor generators}
Q_{jl}= z_j z_l^* ,
\end{equation}
via the generators $\{z_j$, $z^*_j\}$ of $\A(S^7_\theta)$. Indeed,
from equation \eqref{eqn twistor generators} one could infer the
relations on the generators of $\A(\CP^3_\theta)$ from those on the
generators of $\A(S^7_\theta)$. By its very definition
$\A(\CP^3_\theta)$ is the invariant subalgebra of $\A(S^7_\theta)$
under this $\U(1)$-action and equation (\ref{eqn twistor
generators}) defines an inclusion of algebras $$\A(\CP^3_\theta)
\hookrightarrow \A(S^7_\theta),$$ giving a noncommutative principal
bundle with structure group $\U(1)$. We thus have algebra inclusions
\begin{equation}\label{eqn twistor fibration}\A(S^4_\theta) \hookrightarrow \A(\CP^3_\theta)
\hookrightarrow \A(S^7_\theta),\end{equation} with the left-hand
arrow still to be determined. As in the classical case, this
inclusion is not a principal fibration (the `typical fibre' is a copy
of the undeformed $\CP^1$) but we may nevertheless express the generators of
$\A(\CP^3_\theta)$ in terms of the generators of $\A(S^4_\theta)$.
For this we need the non-degenerate map on $\A(\C^4_\theta)$ given
on generators by
\begin{equation}\label{eqn J}J(z_1,z_2,z_3,z_4):= (-z_2^*,z_1^*,-z_4^*,z_3^*)
\end{equation} and extended
as an anti-algebra map. Classically, in doing so we would be identifying the set of
quaternions $\HH$ with the set of $2\times 2$ matrices over $\C$ of
the form $$c_1+c_2j\in \HH ~~\mapsto~~ \begin{pmatrix} c_1&-\bar
c_2\\c_2&\bar c_1 \end{pmatrix}\in \M_2(\C),
$$ and the map $J$
corresponds to right multiplication by the quaternion $j$. In the
deformed case, this very same identification defines the algebra
$\A(\HH^2_\theta)$ to be equal to the algebra $\A(\C^4_\theta)$
equipped with the map $J$ \cite{lprs:ncfi}.

Using the identification of generators (\ref{eqn twistor
generators}) the map $J$ extends to an automorphism of
$\A(\CP^3_\theta)$, given in terms of the matrix generators in
equation (\ref{eqn twistor matrix}) by
\begin{align*}
&J(t_1)  = t_2,  \qquad J(t_2) = t_1,   \qquad J(t_3) = t_4,  \qquad J(t_4) =t_3, \\
&J(x_1)  = -x_1,  \quad J(y_1) =-y_1,  \quad J(x_1^*) =-x_1^*, \quad J(y_1^*) =-y_1^*, \\
&J(x_2)=\mu\, y_2^*,  \quad J(x_3) = -y_3^*, \quad J(x_2^*)=\bar\mu\, y_2, \quad J(x_3^*)=-y_3, \\
&J(y_2)=\bar\mu\, x_2^*,  \quad J(y_3) =-x_3^*, \quad J(y_2^*) = \mu\, x_2, \quad J(y_3^*)=-x_3,
\end{align*}
as required for $J$ to respect the algebra
relations of $\A(\CP^3_\theta)$. The subalgebra fixed by the map $J$ is
precisely $\A(S^4_\theta)$; in fact one has an algebra inclusion 
$\A(S^4_\theta) \hookrightarrow \A(\CP^3_\theta)$ given on
generators by
\begin{equation} \label{eqn twistor to sphere} x \mapsto 2(t_1+t_2-1), \quad
\alpha \mapsto 2(x_2 + \mu\, y_2^*),\quad \beta
\mapsto 2(-x_3^*+y_3),
\end{equation}
with $\mu = \sqrt{\lambda}=\exp{(\pi \ii \theta)}$. In the notation of 
equation (\ref{eqn partial isometry}) we have $Q=|\psi_1\ra\la\psi_1|$,
and we note also that $|\psi_2\ra=|J\psi_1\ra$, so that equation
(\ref{eqn twistor to sphere}) is just the statement that
$$\qp=|\psi_1\ra\la\psi_1|+|\psi_2\ra\la\psi_2|=|\psi_1\ra\la\psi_1|+|J\psi_1\ra\la J\psi_1|=Q+J(Q).$$
This gives us the promised algebraic description of the twistor
fibration (\ref{eqn twistor fibration}): the generators of
$\A(S^4_\theta)$ are identified with the degree one elements of
$\A(\CP^3_\theta)$ of the form $Z+J(Z)$.

\section{The Quantum Conformal Group}
\label{section noncommutative conformal transformations}
Next,  we briefly review the construction of the quantum groups
which describe the symmetries of the spheres $S^4_\theta$ and
$S^7_\theta$ (and the symmetries of the Hopf fibration defined in
Sect.~\ref{section noncommutative hopf fibration}).

\subsection{The quantum groups $\SL_\theta(2,\HH)$ and
$\Sp_\theta(2)$} To begin, we need a noncommutative analogue of the
set of all linear transformations of the quaternionic vector space
$\HH^2_\theta$ defined above. To this end, we define a transformation bialgebra
for the algebra $\A(\HH^2_\theta)$ to be a bialgebra $\mathcal{B}$
such that there is a $*$-algebra map
$\Delta_L:\A(\C^4_\theta)\to \mathcal{B}\otimes
\A(\C^4_\theta)$ commuting with the map $J$ of equation (\ref{eqn
J}). The set of all transformation bialgebras for $\A(\HH^2_\theta)$
forms a category in the natural way; we define the bialgebra
$\A(\M_\theta(2,\HH))$ as the universal initial object in the
category, meaning that whenever $\mathcal{B}$ is a transformation
bialgebra for $\A(\HH^2_\theta)$ there is a morphism of
transformation bialgebras $\A(\M_\theta(2,\HH))\to
\mathcal{B}$ \cite{lprs:ncfi}.
Using the universality property, one finds that
$\A(\M_\theta(2,\HH))$ is the associative algebra generated by the
entries of the following $4\times 4$ matrix:
\begin{equation} \label{eqn defining
M(H)}A=\begin{pmatrix}a_{ij}&b_{ij}\\c_{ij}&d_{ij}\end{pmatrix}=\begin{pmatrix}a_1&-a_2^*&b_1&-b_2^*\\a_2&a_1^*&b_2&b_1^*\\c_1&-c_2^*&d_1&-d_2^*\\c_2&c_1^*&d_2&d_1^*\end{pmatrix}.\end{equation}
With our earlier notation, we think of this matrix as
generated by four quaternion-valued functions, writing
$$a=(a_{ij})=\begin{pmatrix}a_1&-a_2^*\\a_2&a_1^*\end{pmatrix}$$ and
similarly for the other entries $b,c,d$. The coalgebra structure
on $\A(\M_\theta(2,\HH))$ is given by
$$\Delta(A_{ij})=\sum\nolimits_l A_{il}\otimes A_{lj}, \qquad
\ep(A_{ij})=\delta_{ij}$$ for $i,j=1,\ldots,4$, and its
$*$-structure is evident from the matrix (\ref{eqn defining M(H)}).
The coaction $\Delta_L$ is determined to be
\begin{equation} \label{eqn conformal coaction}
\Delta_L:\A(\C^4_\theta)\to
\A(\M_\theta(2,\HH))\otimes\A(\C^4_\theta),\qquad
\Delta_L(\Psi_{ia})=\sum\nolimits_j A_{ij}\otimes
\Psi_{ja},\end{equation} where $\Psi$ is the matrix in equation
(\ref{eqn partial isometry}) (although here we do not assume the
sphere relation and instead think of the entries of $\Psi$ as
generators of the algebra $\A(\C^4_\theta)$). The relations between
the generators of $\A(\M_\theta(2,\HH))$ are found from the
requirement that $\Delta_L$ make $\A(\C^4_\theta)$ into an
$\A(\M_\theta(2,\HH))$-comodule algebra. One computes
\begin{equation}\label{eqn quantum group
relations} \Delta_L(\Psi_{ia}\Psi_{jb})=\sum\nolimits_{km}
(A_{im}A_{jl}-\eta_{ij}\eta_{lm}A_{jl}A_{im})\otimes
\Psi_{ma}\Psi_{lb}\end{equation} and, since the products
$\Psi_{ma}\Psi_{lb}$ may be taken to be all independent as $k,l,a,b$
vary, we must have that
\begin{equation} \label{eqn twisted quantum group
relations}A_{im}A_{jl}=\eta_{ij}\eta_{lm}A_{jl}A_{im}\end{equation}
for $i,j,l,m=1,\ldots,4$. It is not difficult to see that the
algebra generated by the $a_{ij}$ is commutative, as are the
algebras generated by the $b_{ij}$, $c_{ij}$, $d_{ij}$, although
overall the algebra is noncommutative due to some non-trivial
relations among components in different blocks.

Of course, $\A(\M_\theta(2,\HH))$ is not quite a Hopf algebra since
it does not have an antipode. We obtain a Hopf algebra by passing to
the quotient of $\A(\M_\theta(2,\HH))$ by the Hopf $*$-ideal
generated by the element $D-1$, where $D=\det A$ is the formal
determinant of the matrix $A$ in (\ref{eqn defining M(H)}). We
denote the quotient by $\A(\SL_\theta(2,\HH))$, the coordinate
algebra on the quantum group $\SL_\theta(2,\HH)$ of matrices in
$\M_\theta(2,\HH)$ with determinant one, and continue to write the
generators of the quotient as $A_{ij}$. The algebra
$\A(\SL_\theta(2,\HH))$ inherits a $*$-bialgebra structure from that
of $\A(\M_\theta(2,\HH))$ and we use the
determinant to define an antipode $S:\A(\SL_\theta(2,\HH))\to
\A(\SL_\theta(2,\HH))$ as in \cite{lprs:ncfi}. The datum $(\A(\SL(2,\HH)),\Delta,\ep,S)$
constitutes a Hopf $*$-algebra.

The Hopf algebra $\A(\Sp_\theta(2))$ is the
quotient of $\A(\SL_\theta(2,\HH))$ by the two-sided $*$-Hopf ideal
generated by
$$\sum\nolimits_l (A^*)_{li}A_{lj} - \delta_{ij},
\qquad i,j=1,\ldots,4.
$$ In this algebra we have the relations
$A^*A=AA^*=1$, or equivalently $S(A)=A^*$. This Hopf algebra is
the coordinate algebra on the quantum group $\Sp_\theta(2)$, the
subgroup of $\SL_\theta(2,\HH)$ of unitary matrices.

Finally there is an inclusion of algebras
$\A(S^7_\theta) \hookrightarrow \A(\Sp_\theta(2))$ given on
generators by the $*$-algebra map \begin{equation} \label{eqn seven
sphere isomorphism} z_1 \mapsto a_1, \quad z_2 \mapsto a_2, \quad
z_3 \mapsto c_1, \quad z_4 \mapsto c_2.\end{equation} This means
that we may identify the first two columns of the matrix $A$ with
the matrix $\Psi$ of equation (\ref{eqn partial isometry}).
Similarly there is an algebra inclusion $\A(S^4_\theta)
\hookrightarrow \A(\Sp_\theta(2))$ given by
\begin{equation}\label{eqn four sphere isomorphism}x \mapsto a_1a_1^*-a_2a_2^*+c_1c_1^*
-c_2c_2^*, \quad \alpha \mapsto a_1c_1^*-a_2^*c_2, \quad \beta
\mapsto -a_1^*c_2+a_2c_1^*.\end{equation} These inclusions yield
algebra isomorphisms of $\A(S^7_\theta)$ and $\A(S^4_\theta)$ with
certain subalgebras of $\A(\Sp_\theta(2))$ of coinvariants under
coactions by appropriate sub-Hopf algebras, thus realising the
noncommutative spheres as quantum homogeneous spaces for
$\Sp_\theta(2)$. We refer to \cite{lprs:ncfi} for details of these
constructions.

\subsection{Quantum conformal transformations}\label{section quantum conformal transformations}We now review how the quantum groups obtained in the previous
section (co)act on the spheres $S^7_\theta$ and $S^4_\theta$ as `quantum
symmetries'. The coaction
\begin{equation} \Delta_L:\A(\C^4_\theta)\to
\A(\SL_\theta(2,\HH))\otimes\A(\C^4_\theta),\qquad
\Delta_L(\Psi_{ia})=\sum\nolimits_j A_{ij}\otimes \Psi_{ja},\end{equation} is by
construction a $*$-algebra map and so, if we assume that the
quantity
$$r^2:=\sum\nolimits_j z^*_j z_j$$ is invertible with inverse $r^{-2}$, then we may also define an inverse for the quantity
$$\rho^2:=\Delta_L\left(\sum\nolimits_j z^*_j z_j\right)$$ by
$\rho^{-2}:=\Delta_L(r^{-2})$. Inverting $r^2$ corresponds to
deleting the origin in $\C^4_\theta$ and we define the coordinate
algebra of the corresponding subset of $\C^4_\theta$ by
$$\A_0(\C^4_\theta):=\A(\C^4_\theta)[r^{-2}],$$ the
algebra $\A(\C^4_\theta)$ with $r^{-2}$ adjoined. Extending
$\Delta_L$ as a $*$-algebra map gives a well-defined coaction
$$\Delta_L:\A_0(\C^4_\theta)\to \A(\SL_\theta(2,\HH))\otimes
\A_0(\C^4_\theta)$$ for which $\A_0(\C^4_\theta)$ is an
$\A(\SL_\theta(2,\HH))$-comodule algebra.                     %Achange at CMI: line break

Writing $\A_0(\widetilde\C^4_\theta):=\Delta_L(\A_0(\C^4_\theta))$
for the image of $\A_0(\C^4_\theta)$ under $\Delta_L$, both
$\rho^{2}$ and $\rho^{-2}$ are central in the algebra
$\A_0(\widetilde \C^4_\theta)$, since $r^{2}$ and $r^{-2}$ are
central in $\A_0(\C^4_\theta)$.

Now the coaction $\Delta_L$ descends to a coaction of the Hopf
algebra $\A(\Sp_\theta(2))$,
\begin{equation}\label{eqn
unitary coaction}\Delta_L:\A_0(\C^4_\theta)\to \A(\Sp_\theta(2))
\otimes \A_0(\C^4_\theta),\end{equation} by the same formula
(\ref{eqn conformal coaction}) now viewed for the quotient
$\A(\Sp_\theta(2))$. In particular, for this coaction one has
$$(\Psi^*\Psi)_{ab}\mapsto\sum\nolimits_{ijl} (A^*)_{li}A_{ij}\otimes
(\Psi^*)_{al}\Psi_{jb}=\sum\nolimits_{jl}\delta_{lj}\otimes
(\Psi^*)_{al}\Psi_{jb}=1\otimes (\Psi^*\Psi)_{ab},$$ since the generators $A_{ij}$ satisfy the
relations $\sum_i(A^*)_{li}A_{ij}=\delta_{lj}$ in the
algebra $\A(\Sp_\theta(2))$. Then
both $\A(S^4_\theta)$ and $\A(S^7_\theta)$ are
$\A(\Sp_\theta(2))$-comodule algebras, since this coaction preserves
the sphere relations (\ref{eqn foursphere relation}) and (\ref{eqn
sevensphere relation}).

In contrast, the spheres $S^7_\theta$ and $S^4_\theta$ are not preserved
under the coaction of the larger quantum group $\SL_\theta(2,\HH)$.
Although defined on the algebra $\A_0(\C^4_\theta)$, the coaction
$\Delta_L$ of $\A(\SL_\theta(2,\HH))$ is not well-defined on the
seven-sphere $\A(S^7_\theta)$ since it does not preserve the sphere
relation $r^2=1$ of equation (\ref{eqn sevensphere relation}). By
definition, we have instead that $\Delta_L(r^2)=\rho^2$, meaning
that the coaction of $\A(\SL_\theta(2,\HH))$ `inflates' the sphere
$\A(S^7_\theta)$\ \cite{lprs:ncfi}.
Since $r^2$ is a central element of $\A_0(\C^4_\theta)$, we may
evaluate it as a positive real number. The result is the coordinate
algebra of a noncommutative sphere $S^7_{\theta,r}$ of radius $r$;
as this radius varies in $\A_0(\C^4_\theta)$, it sweeps out a family
of seven-spheres. Similarly, evaluation of the central element
$\rho^2$ in $\A_0(\widetilde \C^4_\theta)$ yields the coordinate algebra
of a noncommutative sphere $\widetilde S^7_{\theta,\rho}$ of radius
$\rho$ and, as the value of $\rho$ varies in $\A_0(\widetilde\C^4_\theta)$, it sweeps out another family of seven-spheres. The
coaction $\Delta_L$ of $\A(\SL_\theta(2,\HH))$ on
$\A_0(\C^4_\theta)$ serves to map the family parameterised by $r^2$
onto the family parameterised by $\rho^2$.

A similar fact is found for the generators $\alpha$, $\beta$, $x$ of
the four-sphere algebra $\A(S^4_\theta)$. The coaction of
$\A(\SL_\theta(2,\HH))$ does not preserve the sphere relation but
gives instead that
$$\Delta_L(\alpha^*\alpha+\beta^*\beta+x^2)=\rho^4,$$ and the four-sphere $S^4_\theta$ is
also inflated. Let us write $\A(\mathcal{Q}_\theta)$ for the subalgebra of
$\A_0(\C^4_\theta)$ generated by $\alpha$, $\beta$, $x$ and their
conjugates. Then as $r^4$ varies in $\A(\mathcal{Q}_\theta)$, we get a family
of noncommutative four-spheres. Similarly, we define $\tilde
\alpha:=\Delta_L(\alpha)$, $\tilde \beta:=\Delta_L(\beta)$, $\tilde
x:=\Delta_L(x)$ and so forth, and write $\A(\widetilde{\mathcal{Q}}_\theta)$
for the subalgebra of $\A_0(\widetilde\C^4_\theta)$ that they
generate. It is precisely the $\SU(2)$-invariant subalgebra of
$\A_0(\widetilde\C^4_\theta)$, and as $\rho^4$ varies in
$\A(\widetilde{\mathcal{Q}}_\theta)$ we get another family of noncommutative
four-spheres. The coaction of the quantum group
$\A(\SL_\theta(2,\HH))$ maps the family parameterised by
$r^4$ onto the family parameterised by $\rho^4$.

Thus there is a family of $\SU(2)$-principal fibrations
given by the algebra inclusion $\A(\mathcal{Q}_\theta)\hookrightarrow
\A_0(\C^4_\theta)$, the family being parameterised by the function
$r^2$. For a fixed value of $r^2$ we get an $\SU(2)$ principal
bundle $S^7_{\theta,r} \to S^4_{\theta,r^2}$. Similarly,
the algebra inclusion $\A(\widetilde{\mathcal{Q}}_\theta) \hookrightarrow
\A_0(\widetilde \C^4_\theta)$ defines a family of $\SU(2)$-principal
fibrations parameterised by the function $\rho^2$. The above
construction shows that the coaction of the quantum group
$\A(\SL_\theta(2,\HH))$ carries the former family of principal fibrations
onto the latter.

All of this means that, as things stand, we cannot use the
presentations of $\A(S^4_\theta)$ and $\A(S^7_\theta)$ of
Sect.~\ref{section noncommutative hopf fibration} to give a
well-defined coaction of $\A(\SL_\theta(2,\HH))$, since the sphere
relations we use to define them are not preserved by the coaction.
Rather we should work with the families of spheres all at once (this
is the price we have to pay for working with the coaction of a Hopf
algebra rather than the action of a group). To do this, we note that
the algebra $\A(S^4_\theta)$ may be identified with the subalgebra
of $\A_0(\C^4_\theta)$ generated by $r^{-2}\alpha,$ $r^{-2}\beta$,
$r^{-2}x$, together with their conjugates, since the sphere relation
\begin{equation}\label{eqn modified sphere relation}(r^{-2}\alpha)(r^{-2}\alpha)^*+(r^{-2}\beta)(r^{-2}\beta)^* +
(r^{-2}x)^2=1\end{equation} is automatically satisfied in
$\A_0(\C^4_\theta)$. The result of doing so is that we have a
well-defined coaction,
$$\Delta_L:\A(S^4_\theta) \to \A(\SL_\theta(2,\HH))\otimes
\A(S^4_\theta), $$ defined on the generators $r^{-2}\alpha$,
$r^{-2}\beta$, $r^{-2}x$ and their conjugates, with the sphere
relation (\ref{eqn modified sphere relation}) now preserved by
$\Delta_L$. In this way, we think of $\SL_\theta(2,\HH))$ as the
quantum group of conformal transformations of $S^4_\theta$.

In these new terms, the construction of the defining projector for
$\A(S^4_\theta)$ needs to be modified only slightly. We now take the
normalised matrix
\begin{equation}\label{eqn modified partial isometry}
\Psi=r^{-1}\begin{pmatrix}
z_1&z_2&z_3&z_4\\-z_2^*&z_1^*&-z_4^*&z_3^*\end{pmatrix}^{\textup{t}},
\end{equation} at the price of including the generator $r^{-1}$ as
well (not a problem in the smooth closure \cite{lprs:ncfi}). Thanks
to the relation (\ref{eqn modified sphere relation}), we still have
$\Psi^*\Psi=1$ and the required projector is
\begin{equation} \label{eqn modified basic instanton projector}
\qp :=\Psi\Psi^*=\tfrac{1}{2}r^{-2}\begin{pmatrix} r^2+x & 0 & \alpha
& -\bar \mu \,\beta^* \\ 0 & r^2+x & \beta & \mu \, \alpha^* \\
\alpha^* & \beta^* & r^2-x & 0 \\ -\mu \,\beta & \bar\mu \, \alpha &
0 & r^2-x
\end{pmatrix}.\end{equation}
By the above discussion, the coaction
$\Delta_L$ of $\A(\SL_\theta(2,\HH))$ is now well-defined on the
algebra generated by the entries of this matrix. Writing $\widetilde\Psi_{ia}:=\Delta_L(\Psi_{ia})$,
the image of $\qp$ under $\Delta_L$
is computed to be
\begin{equation} \label{eqn coacted instanton projector}
\tilde{\qp}:=\widetilde\Psi\widetilde\Psi^*=\tfrac{1}{2}\rho^{-2}\begin{pmatrix}
\rho^2+\tilde x & 0 &
\tilde \alpha & -\bar \mu \,\tilde \beta^* \\ 0 & \rho^2+x & \tilde \beta & \mu \, \tilde \alpha^* \\
\tilde \alpha^* & \tilde \beta^* & \rho^2-x & 0 \\ -\mu \,\tilde \beta
& \bar\mu \, \tilde \alpha & 0 & \rho^2-\tilde x
\end{pmatrix}.\end{equation}
The entries of these projectors
generate respectively subalgebras of $\A_0(\C^4_\theta)$ and
$\A_0(\widetilde \C^4_\theta)$, each parameterising the families of
noncommutative four-spheres discussed above.

Finally, we observe that similar statements may be made about the
$\U(1)$-principal fibration $S^7_\theta \to \CP^3_\theta$. We do not
need a sphere relation in order to define the coordinate algebra
$\A(\CP^3_\theta)$: in Sect.~\ref{section noncommutative twistor
space} it was merely convenient to do so. Instead, we may identify
$\A(\CP^3_\theta)$ as the $\U(1)$-invariant subalgebra of
$\A_0(\C^4_\theta)$ generated by elements $t_1=r^{-2}z_1z_2^*$,
$x_1=r^{-2}z_1z_2^*$, $x_2=r^{-2}z_1z_3^*$, $x_3= r^{-2}z_1z_4^*$
and so forth.

\section{A Noncommutative ADHM construction} \label{section adhm construction}
There is a well-known solution to the problem of constructing
instantons on the classical four-sphere $S^4$ which goes under the
name of ADHM construction. Techniques of linear algebra are used to
construct vector bundles over twistor space $\CP^3$, which are in
turn put together to construct a vector bundle over $S^4$ equipped
with an instanton connection. It is known that all such
connections are obtained in this way \cite{adhm:ci,ma:gymf}.

Our goal here is to generalise the ADHM method to a deformed version
which constructs instantons on the noncommutative sphere
$S^4_\theta$. The classical construction may be obtained from our
deformed version by setting $\theta=0$. As usual our approach stems
from writing the classical construction in a dualised language which
does not depend on the commutativity of the available function
algebras, although here the situation is not as straightforward as
one might first expect. The deformed construction is rather more
subtle than it is in the commutative case and produces
noncommutative `families' of instantons.

\subsection{A noncommutative space of monads}
%In this section we shall need the following notation.
The algebra $\A(\C^4_\theta)$ has a natural $\ZZ$-grading given by
assigning to its generators the degrees $$\textup{deg}(z_j)=1, \quad
\textup{deg}(z_j^*)=-1, \qquad j=1,\ldots,4,$$ which results in a
decomposition $\A(\C^4_\theta)=\oplus_{n \in \ZZ} \A_n$. Then for
each $r \in \ZZ$ there is a `degree shift' map from
$\A(\C^4_\theta)$ to itself whose image we denote
$\A(\C^4_\theta)(r)$; by definition the degree $n$ component of
$\A(\C^4_\theta)(r)$ is  $\A_{r+n}$.

Similarly, if a given $\A(\C^4_\theta)$-module $\E$ is $\ZZ$-graded, we
denote the degree-shifted modules  by $\E(r)$, $r \in \ZZ$. In
particular, for each finite dimensional vector space $H$ the
corresponding free right module $H \otimes \A(\C^4_\theta)$ is
$\ZZ$-graded by the grading on $\A(\C^4_\theta)$, and the shift maps
on $\A(\C^4_\theta)$ induce the shift maps on $H \otimes
\A(\C^4_\theta)$.

The input data for the classical ADHM construction of $\SU(2)$
instantons with topological charge $k$ is a {\em monad}, by which we
mean a sequence of free right modules over the algebra $\A(\C^4)$,
\begin{equation} \label{eqn module monad} H\otimes \A(\C^4)(-1) \xrightarrow{\sigma_z}
K\otimes \A(\C^4) \xrightarrow{\tau_z} L\otimes
\A(\C^4)(1),\end{equation} where $H$, $K$ and $L$ are complex vector
spaces of dimensions $k$, $2k+2$ and $k$ respectively. The arrows
$\sigma_z$ and $\tau_z$ are $\A(\C^4)$-module homomorphisms assumed
to be such that $\sigma_z$ is injective, $\tau_z$ is surjective and
that the composition $\tau_z \sigma_z=0$. This is the usual approach
in algebraic geometry \cite{oss:vb}, although here we work with
$\A(\C^4)$-modules, {\em i.e.} global sections of vector bundles,
rather than with locally-free sheaves.

The degree shifts signify we think of $\sigma_z$ and $\tau_z$
respectively as elements of $H^*\otimes K \otimes \A_1$ and
$K^*\otimes L\otimes \A_1$, where $\A_1$ is the degree one component
of $\A(\C^4)$ (the vector space spanned by the generators
$z_1,\ldots,z_4$). This means that alternatively we may think of
them as linear maps \begin{equation}\label{eqn param classical
maps}\sigma_z:H\times \C^4 \to K, \qquad \tau_z:K \times \C^4 \to
L,\end{equation} thus recovering the more explicit geometric
approach of \cite{adhm:ci}.

Our goal in this section is to give a description of a monad of the
form (\ref{eqn module monad}) in an algebraic framework which allows
the possibility of the algebra $\A(\C^4_\theta)$ being
noncommutative. In this setting, we require the maps $\sigma_z$ and
$\tau_z$ to be parameterised by the noncommutative space
$\C^4_\theta$ rather than by the classical space $\C^4$, as was the
case in equation (\ref{eqn param classical maps}). Our first task
then is to find an analogue of the space of linear module maps
$H\otimes \A(\C^4_\theta)(-1) \to K \otimes \A(\C^4_\theta)$.

Following a general strategy \cite{wa:qfm}, we define
$\A(\tmM_\theta(H,K))$ to be the universal algebra for which there
is a morphism of right $\A(\C^4_\theta)$-modules,
$$\sigma_z: H\otimes \A(\C^4_\theta)(-1) \to
\A(\tmM_\theta(H,K))\otimes K \otimes \A(\C^4_\theta),$$ which is
linear in the generators $z_1,\ldots,z_4$ of $\A(\C^4_\theta)$. By
this  we mean that whenever $\mathcal{B}$ is an algebra satisfying
these properties there exists a morphism of algebras
$\phi:\A(\tmM_\theta(H,K))\to\mathcal{B}$ and a
commutative diagram $$\begin{CD} H\otimes \A(\C^4_\theta)(-1)
 @>\sigma_z>> \A(\tmM_\theta(H,K))\otimes K \otimes \A(\C^4_\theta)
\\ @VV\id V @VV\phi\otimes \id V \\ H\otimes \A(\C^4_\theta)(-1) @>\sigma_z'>> \mathcal{B}\otimes K \otimes
\A(\C^4_\theta)\end{CD}$$ of right $\A(\C^4_\theta)$-modules, with
$\sigma_z'$ denoting the corresponding map for the algebra
$\mathcal{B}$.

Choosing a basis $(u_1,\ldots,u_k)$ for the vector space $H$ and a
basis $(v_1,\ldots,v_{2k+2})$ for the vector space $K$, the algebra $\A(\tmM_\theta(H,K))$
is generated by the matrix elements
$$\{M_{ab}^\alpha~|~a=1,\ldots,2k+2,~b=1,\ldots,k,~\alpha=1,\ldots,4\},$$
which define a map $\sigma_z$, expressed on simple tensors by
\begin{equation}\label{eqn nc module map}
\sigma_z: u_b\otimes Z \mapsto \sum\nolimits_{a,\alpha}
M_{ab}^\alpha \otimes v_a \otimes z_\alpha Z, \qquad Z\in
\A(\C^4_\theta).
\end{equation} In
more compact notation, for each $\alpha$ we arrange these elements
into a $(2k+2)\times k$ matrix $M^\alpha=(M^\alpha_{ab})$, so that
with respect to the above bases, $\sigma_z$ may be written
\begin{equation} \label{eqn classical sigma}
\sigma_z=\sum\nolimits_\alpha M^\alpha \otimes z_\alpha .
\end{equation}

To find the relations in the algebra $\A(\tmM_\theta(H,K))$, let us
write $(\hat u_1,\ldots,\hat u_k)$ for the basis of $H^*$ which is
dual to $(u_1,\ldots,u_k)$ and write $(\hat v_1,\ldots,\hat
v_{2k+2})$ for the basis of $K^*$ dual to $(v_1,\ldots,v_{2k+2})$.
Then the map (\ref{eqn nc module map}) has an equivalent dual
description (also denoted $\sigma_z$) in terms of the dual vector
spaces $H^*$, $K^*$ as
\begin{equation}\label{eqn dual nc module map}
\sigma_z:\hat v_a\otimes Z \mapsto \sum\nolimits_{b,\alpha}
M^\alpha_{ab}\otimes \hat u_b \otimes z_\alpha Z,
\end{equation}
and extended as an $\A(\C^4_\theta)$-module map. The functionals
$\hat u_b$, $\hat v_a$ together with their conjugates $\hat u_b^*$, $\hat
v_a^*$ generate the coordinate algebras of $H$ and $K$ respectively.
It is only
natural to require that (\ref{eqn dual nc module map}) be an algebra
map.

\begin{prop}\label{prop universal matrix relations}
With $(\eta_{\alpha \beta})$ the
matrix (\ref{eqn eta matrix}) of deformation parameters, the matrix
elements $M^\alpha_{ab}$ enjoy the relations
\begin{equation}\label{eqn sigma relations}M^\alpha_{ab}M^{\beta}_{cd}=\eta_{\beta
\alpha}M^\beta_{cd}M^{\alpha}_{ab}\end{equation} for each
$a,c=1,\ldots,2k+2$, each $b,d=1,\ldots,k$ and each $\alpha,
\beta=1,\ldots,4$.
\end{prop}
\proof The requirement that (\ref{eqn dual nc module map}) is an
algebra map means that in degree one we need $\sigma_z(\hat v_a\hat
v_c)=\sigma_z(\hat v_c\hat v_a)$ for all $a,c=1,\ldots,2k+2$, which
translates into the statement that
$$
\sum\nolimits_{b,d,\alpha,\beta} M_{ab}^\alpha M_{cd}^\beta \otimes
\hat u_b \hat u_d \otimes z_\alpha
z_\beta=\sum\nolimits_{b,d,\alpha,\beta} M_{cd}^\beta M_{ab}^\alpha
\otimes \hat u_d \hat u_b \otimes z_\beta z_\alpha
$$
for all $a,c=1,\ldots,2k+2$. Using in turn the relations (\ref{eqn
nc params}) and the fact that the generators $\hat u_b$, $\hat u_d$
commute for all values of $b,d$, this equation may be rearranged to
give $$\sum\nolimits_{b,d,\alpha,\beta} \left( M_{ab}^\alpha
M_{cd}^\beta -\eta_{\beta\alpha}M_{cd}^\beta M_{ab}^\alpha\right)
\otimes \hat u_b \hat u_d \otimes z_\alpha z_\beta=0.$$ Since for $b
\leq d$ and $\alpha \leq \beta$ the quantities $\hat u_b \hat u_d
\otimes z_\alpha z_\beta$ may all be taken to be independent, we
must have that their coefficients are all zero, leading to the
stated relations.\endproof

The above proposition simply says that the entries of a given matrix
$M^\alpha$ all commute, whereas the relations between the entries of
the matrices $M^\alpha$ and $M^\beta$ are determined by the
deformation parameter $\eta_{\beta \alpha}$. Hence the algebra
$\A(\tmM_\theta(H,K))$ is generated by the $M^\alpha_{ab}$ subject
to the relations (\ref{eqn sigma relations}). The algebra
$\A(\tmM_{\theta=0}(H,K))$ is commutative and parameterises the
space of all possible maps $\sigma_z$, since for each point $x \in
\tmM_{\theta=0}(H,K)$ there is an evaluation map,
$$
\textup{ev}_x:\A(\tmM_{\theta=0}(H,K))\to \C,
$$
which yields an $\A(\C^4)$-module homomorphism
\begin{align*}
(\textup{ev}_x\otimes \textup{id})\sigma_z &:H\otimes
\A(\C^4)(-1)\to K\otimes \A(\C^4), \\ (\textup{ev}_x\otimes
\textup{id})\sigma_z &:=\sum \textup{ev}_x(M^\alpha_{ab})\otimes
z_\alpha.
\end{align*}
When $\theta$ is different from zero, there need not be enough
evaluation maps available. Nevertheless, we think of
$\A(\tmM_\theta(H,K))$ as a noncommutative family of maps
parameterised by the noncommutative space $\tmM_\theta(H,K)$.

\begin{rem}
\textup{Since we constructed $\A(\tmM_\theta(H,K))$ through the
minimal requirement that $\sigma_z$ is an algebra map, it is indeed
the universal algebra with the required properties. This means that
our interpretation of $\A(\tmM_\theta(H,K))$ as a noncommutative
family of maps is in agreement with the approaches of
\cite{wa:qfm,sw:proc,sol:qfm} for quantum families of maps
parameterised by noncommutative spaces. Moreover, it also agrees
with the definition of algebras of rectangular quantum matrices
discussed in \cite{mm:qumat}. It may also be viewed as a kind of
`comeasuring' as introduced in \cite{maj:bdg}, but now for modules
instead of algebras.}\end{rem}

Thus we have a noncommutative analogue of the space of all maps
$\sigma_z$. A similar construction works for the maps
$\tau_z$: there is a universal algebra $\A(\tmM_\theta(K,L))$
generated by matrix elements $N^\alpha_{ba}$ for labels $b=1,\ldots,k$,
$a=1,\ldots,2k+2$ and $\alpha=1,\ldots,4$, here coming from a map
\begin{equation}\label{eqn nc tau expression}\tau_z:v_a\otimes Z\mapsto
\sum\nolimits_{b,\alpha} N_{ba}^\alpha \otimes w_b\otimes z_\alpha
Z,\end{equation} having chosen a basis $(w_1,\ldots,w_k)$ for the
vector space $L$. Dually, the requirement that $\tau_z$ be an
algebra map from the coordinate algebra of $L$ to the coordinate
algebra of $K$ results in relations for the generators of the
algebra $\A(\tmM_\theta(K,L))$,
\begin{equation}\label{eqn-tau-relations}
N^\alpha_{ba}N^\beta_{dc}=\eta_{\beta\alpha}N^\beta_{dc}N^\alpha_{ba},
\end{equation}
which are the parallel of conditions \eqref{eqn sigma relations} for the algebra
$\A(\tmM_\theta(H,K))$.

To complete the monad picture we finally require that the
composition of the maps $\sigma_z$ and $\tau_z$ be zero. In the
dualised format the composition is easily dealt with as the
composition as a map from the coordinate algebra of $L$ to that of
$H$, with the product appearing as part of a general procedure for
`gluing' quantum matrices \cite{mm:qumat}. By this we mean that the
composition $\vartheta_z:=\tau_z\circ\sigma_z$ is given in terms of
an algebra-valued $k\times k$ matrix, the product of a $k\times
(2k+2)$ matrix with a $(2k+2)\times k$ matrix. Explicitly, the map
is
$$\vartheta_z:H\otimes\A(\C^4_\theta)(-1)\to\A(\widetilde{\mathcal{M}}_\theta(H,L))\otimes L\otimes
\A(\C^4_\theta)(1),$$
$$\vartheta_z: \hat w_a\otimes Z \mapsto \sum\nolimits_{b,\alpha,\beta}
T^{\alpha,\beta}_{ab}\otimes \hat w_b \otimes z_\alpha z_\beta Z,$$
where $\A(\widetilde{\mathcal{M}}_\theta(H,L))$ is the coordinate
algebra generated by the matrix elements $T^{\alpha,\beta}_{ab}$ for
$\alpha,\beta=1,\ldots,4$ and $a,b=1,\ldots,k$.
The matrix
multiplication $(\tau_z,\sigma_z)\mapsto \vartheta_z$ now appears as
a `coproduct'
$$\A(\widetilde{\mathcal{M}}_\theta(H,L))\to\A(\widetilde{\mathcal{M}}_\theta(K,L))\otimes\A(\widetilde{\mathcal{M}}_\theta(H,K)),$$
$$T^{\alpha,\beta}_{cd}:=\sum\nolimits_b
N^{\alpha}_{cb}\otimes M^{\beta}_{bd},\qquad \alpha,\beta=1,\ldots,4,~c,d=1,\ldots,k.$$
The condition $\tau_z\sigma_z=0$ is thus that the image of
this map in
$\A(\widetilde{\mathcal{M}}_\theta(K,L))\otimes\A(\widetilde{\mathcal{M}}_\theta(H,K))$
is zero; this is established by the following proposition.

\begin{prop}\label{prop monad composition relation} The condition $\tau_z\sigma_z=0$ is equivalent to the
requirement that
\begin{equation}\label{eqn monad relations}
\sum\nolimits_r (N^\alpha_{br}M^\beta_{rd}+\eta_{\beta
\alpha}N^\beta_{br}M^\alpha_{rd})=0
\end{equation}
for all $b,d=1,\ldots,k$ and all $\alpha, \beta = 1,\ldots,4$.
\end{prop}

\proof In terms of algebra-valued matrices the map $\tau_z\sigma_z$
is computed as the composition of the duals of the maps \eqref{eqn
nc module map} and \eqref{eqn nc tau expression}, following the
discussion above, to be equal to
$$
(\tau_z\sigma_z)_{bd}=\sum\nolimits_{r,\alpha,\beta} N_{br}^\alpha
M_{rd}^\beta \otimes z_\alpha z_\beta .
$$
Equating to zero the coefficients of the linearly independent
generators $z_\alpha z_\beta$ for $\alpha \leq \beta$ gives the
relations as stated.\endproof

The conditions in equation (\ref{eqn monad relations}) may be
expressed more compactly in terms of products of matrices as
$$N^\alpha M^\beta +\eta_{\beta\alpha}N^\beta M^\alpha=0,$$ for
$\alpha, \beta=1,\ldots,4$ (and as in \eqref{eqn-tau-relations} there is no sum over
$\alpha$ and $\beta$ in this expression).

\begin{defn}\label{defn univ monad algebra}
Define $\A(\tmM_{\theta;k})$ to be the algebra generated by the
matrix elements $M^\alpha_{ab}$ and $N^\beta_{ba}$ subject to the
relations
$$M^\alpha_{ab}M^{\beta}_{cd}=\eta_{\beta
\alpha}M^\beta_{cd}M^{\alpha}_{ab}, \quad
N^\alpha_{ba}N^\beta_{dc}=\eta_{\beta\alpha}N^\beta_{dc}N^\alpha_{ba},
$$
as well as the relations
$$\sum\nolimits_r(N^\alpha_{dr}M^\beta_{rb}+\eta_{\beta\alpha}N^\beta_{br}M^\alpha_{rd})=0
$$
for all $\alpha,\beta=1,\ldots,4$, all $b,d=1,\ldots,k$ and all
$a,c=1,\ldots,2k+2$.\end{defn}

The noncommutative algebra $\A(\tmM_{\theta;k})$ is by construction
universal amongst all algebras having the property that the
resulting maps $\sigma_z$ and $\tau_z$ are algebra maps which
compose to zero. Our interpretation is that for fixed $k$ the
collection of monads over $\C^4_\theta$ is parameterised by the
noncommutative space which is `dual' to this algebra.

\subsection{The subspace of self-dual monads} In the classical
case, the input datum of a monad is by itself insufficient to
construct bundles over the four-sphere $S^4$. To achieve this, one
must incorporate the quaternionic structure afforded by the map $J$
as in \eqref{eqn J} (in the classical limit) and ensure that the
monad is compatible with this extra structure. The same is true in
the noncommutative case, as we shall see presently.

Given the pair of maps constructed in the previous section,
\begin{align*}
\sigma_z:H\otimes \A(\C^4_\theta)(-1) &\to
\A(\tmM_\theta(H,K))\otimes K \otimes \A(\C^4_\theta), \\
\tau_z:K \otimes \A(\C^4_\theta) & \to \A(\tmM_\theta(K,L))\otimes L
\otimes \A(\C^4_\theta)(1),
\end{align*}
we firstly note that the anti-algebra map $J$ in \eqref{eqn J}
induces a new pair of maps,
\begin{align}\label{eqn J applied to monad1}
\sigma_{J(z)}&:H\otimes J\left(\A(\C^4_\theta)(-1)\right) \to
\A\left(\tmM_\theta(H,K)\right)\otimes K \otimes
J\left(\A(\C^4_\theta)\right), \nonumber \\
\sigma_{J(z)}&:=\sum_\alpha M^\alpha \otimes J(z_\alpha),
\end{align}
and
\begin{align}\label{eqn J applied to monad2}
\tau_{J(z)}&:K \otimes J\left(\A(\C^4_\theta)\right) \to
\A\left(\tmM_\theta(K,L)\right)\otimes L
\otimes J\left(\A(\C^4_\theta)(1)\right), \nonumber \\
\tau_{J(z)}&:=\sum_\alpha N^\alpha \otimes J(z_\alpha).
\end{align}
Here, $J\left(\A(\C^4_\theta)\right)$ is the left
$\A(\C^4_\theta)$-module induced by the anti-algebra map $J$ and
$\sigma_{J(z)}$, $\tau_{J(z)}$ are homomorphisms of left
$\A(\C^4_\theta)$-modules. We may also take the adjoints of the
above maps. To make sense of this, we need to add to our picture the
matrix elements $M^\alpha_{ab}{}^*$, so that the adjoint of
$\sigma_z$ is
\begin{equation}\label{eqn adjoint nc module map}
\sigma_z^\star: v_a\otimes Z
\mapsto\sum\nolimits_{b,\alpha}M^\alpha_{ab}{}^* \otimes u_b \otimes
z_\alpha^*Z, \qquad Z \in \A(\C^4_\theta),
\end{equation} where $a=1,\ldots,2k+2$, $b=1,\ldots,k$ and $\alpha=1,\ldots,4$. Let us denote by
$M^\alpha{}^\dag$ the $k\times (2k+2)$ matrix with entries
$(M^\alpha{}^\dag)_{ba}=M^\alpha_{ab}{}^*$. Then with respect to the
above choice of bases, the adjoint map $\sigma_z^\star$ may be
written more compactly as
$$
\sigma_z^\star=\sum\nolimits_\alpha M^\alpha{}^\dag \otimes
z_\alpha^*.
$$ Similarly, we add the matrix
elements $N^\alpha_{dc}{}^*$ and write
$(N^\alpha{}^\dag)_{cd}=N^\alpha_{dc}{}^*$, so that the adjoint of
$\tau_z$ is $$\tau_z^\star: w_b\otimes Z \mapsto
\sum\nolimits_{a,\alpha} N_{ba}^\alpha{}^* \otimes v_a \otimes
z_\alpha^*Z,$$ or $\tau_z^\star=\sum_\alpha N^\alpha{}^\dag \otimes
z_\alpha^*$ in compact notation. The elements $M^\alpha_{ab}{}^*$
are the generators of the algebra $\A(\tmM_\theta(K^*,H^*))$,
whereas the elements $N^\alpha_{dc}{}^*$ are the generators the
algebra $\A(\tmM_\theta(L^*,K^*))$. Applied to equations (\ref{eqn J
applied to monad1}) and (\ref{eqn J applied to monad2}), all of this
yields a pair of homomorphisms of right $\A(\C^4_\theta)$-modules
\begin{align*}
\sigma_{J(z)}^\star:K^*\otimes J\left(\A(\C^4_\theta)\right)^* &\to
\A\left(\tmM_\theta(K^*,H^*)\right)\otimes H^* \otimes
J\left(\A(\C^4_\theta)\right)^*(1),
%\label{eqn J applied to monad3}
\\
\tau_{J(z)}^\star:L^* \otimes J\left(\A(\C^4_\theta)\right)^*(-1) &
\to \A\left(\tmM_\theta(L^*,K^*)\right)\otimes K^* \otimes
J\left(\A(\C^4_\theta)\right)^*,
%\label{eqn J applied to monad4}
\end{align*} defined respectively by
$$\sigma_{J(z)}^\star=\sum_\alpha M^\alpha{}^\dag \otimes
J(z_\alpha)^*,\qquad \tau_{J(z)}^\star=\sum_\alpha N^\alpha{}^\dag
\otimes J(z_\alpha)^*.$$ Of course, we may identify the vector
spaces $H$ and $L^*$ through the basis isomorphism $u_b \mapsto \hat
w_b$ for each $b=1,\ldots,k$. Similarly the isomorphism $v_a \mapsto
\hat v_a$ for $a=1,\ldots,2k+2$ gives an identification of the
vector space $K$ with its dual $K^*$. Also, the right module
$J(\A(\C^4_\theta))^*$ may be identified with $\A(\C^4_\theta)$ by
the composition of the map $J$ with the involution $*$ (noting that
this identification is not the identity map). Through these
identifications, we may think of $\sigma_{J(z)}^\star$ and
$\tau_{J(z)}^\star$ as module homomorphisms
\begin{align*}
\tau_{J(z)}^\star:H \otimes \A(\C^4_\theta)(-1) &\to
\A\left(\tmM_\theta(H,K)\right)\otimes K \otimes\A(\C^4_\theta), \\
%\label{eqn J applied to monad5}
\sigma_{J(z)}^\star:K\otimes \A(\C^4_\theta) &\to
\A\left(\tmM_\theta(K,L)\right)\otimes L \otimes\A(\C^4_\theta)(1).
%\label{eqn J applied to monad6}
\end{align*}
It is straightforward to check that we now have
$\sigma_{J(z)}^\star\tau_{J(z)}^\star=0$ and so all of this means
that the maps $\sigma_{J(z)}^\star$ and $\tau_{J(z)}^\star$ also
give a parameterisation of the noncommutative space of monads,
albeit a different parameterisation from the one we started with. In
the classical case the above procedure applied to a given monad
again yields a monad, although it is not necessarily the one we
started with. If fact, in the classical case, one is interested only
in the subset of monads which are invariant under the above
construction, namely the monad obtained by applying $J$ and
dualising is required to be isomorphic to the one we start with
(this is the sense in which we require monads to be compatible with
$J$). We call such monads {\em self-dual}. In our algebraic
framework, where we work not with specific monads but rather with
the (possibly noncommutative) space $\tmM_{\theta;k}$ of all monads,
this extra requirement is encoded as follows.

\begin{prop}\label{prop matrix monad conditions} The space of self-dual monads is parameterised by the algebra $\A(\tmM^{SD}_{\theta;k})$, the quotient of
the algebra $\A(\tmM_{\theta;k})$ by the further relations
\begin{equation}\label{eqn matrix reality conditions}N^1=-M^2{}^\dag,
\quad N^2=M^1{}^\dag, \quad N^3=-M^4{}^\dag, \quad
N^4=M^3{}^\dag.\end{equation}\end{prop}

\proof The condition that the maps $\sigma_z$ and $\tau_z$ should
parameterise self-dual monads is that $\sigma_z=\tau_{J(z)}^\star$,
equivalently that $\tau_z=-\sigma_{J(z)}^\star$. In terms of the
matrices $M^\alpha$, $N^\alpha$, the former condition reads
\begin{equation}\label{eqn expanded reality condition}
\sum\nolimits_\alpha M^\alpha \otimes z_\alpha=\sum\nolimits_\alpha
N^\alpha{}^\dag \otimes J(z_\alpha)^*.\end{equation} Equating
coefficients of generators of $\A(\C^4_\theta)$ in each of these
equations yields the extra relations as stated. \endproof

\begin{rem}\label{rem orth columns}
\textup{The identification of the vector space $K$ with its dual
$K^*$ means that the module $K \otimes \A(\C^4_\theta)$ acquires a
bilinear form given by
\begin{equation}\label{eqn bilinear form}
(\xi,\eta):=\la J\xi|\eta\ra=\sum\nolimits_a (J\xi)^*_a\eta_a
\end{equation}
for $\xi=(\xi_a)$ and $\eta=(\eta_a) \in K\otimes \A(\C^4_\theta)$,
with $\la \,\cdot\,|\,\cdot \,\ra$ the canonical Hermitian structure
on $K\otimes \A(\C^4_\theta)$. The monad condition, which now reads
$$
0=\tau_z\sigma_z=-\sigma_{J(z)}^\star \sigma_z,
$$  translates into the more practical condition that
the columns of the matrix $\sigma_z$ (equivalently the rows of
$\tau_z$) are orthogonal with respect to the form $(\,\cdot \,,
\,\cdot \,)$. }
\end{rem}
\noindent
Moreover, we see that
\begin{align*}0&=\tau_{z+J(z)}\sigma_{z+J(z)}=\tau_z\sigma_z +
\tau_z\sigma_{J(z)}+\tau_{J(z)}\sigma_z+\tau_{J(z)}\sigma_{J(z)}
=\tau_z\sigma_{J(z)}+\tau_{J(z)}\sigma_z
\\&=-\sigma^\star_{J(z)}\sigma_{J(z)}+\sigma_z^\star
\sigma_z\end{align*} so that in the matrix algebra $\M_k(\C)\otimes
\A(\tmM^{SD}_{\theta;k})\otimes \A(\C^4_\theta)$ we have also
$$
\sigma_{J(z)}^\star\sigma_{J(z)}=\sigma_z^\star\sigma_z.
$$
\begin{rem}\label{rem
number of params} \textup{The above identifications of vector spaces
$H\cong L^*$ and $K\cong K^*$ yield an identification of
$\A(\tmM_\theta (H,K))$ with $\A(\tmM_\theta (L^*,K^*))$ and hence a
reality structure on the generators $M^\alpha_{ab}$. It follows that
the space of self-dual monads is parameterised by a total of
$4k(2k+2)$ generators $M^\alpha_{ab}$. As already remarked, the
condition $\sigma_{J(z)}^\star \sigma_z=0$ is equivalent to
demanding that the columns of $\sigma_z$ are pairwise orthogonal
with respect to the bilinear form $(\,\cdot\,,\,\cdot\,)$ and, since
$\sigma_z$ has $k$ columns, this yields $\frac{1}{2}k(k-1)$ such
orthogonality conditions. Now as in Prop.~\ref{prop monad
composition relation} we may equate to zero the coefficients of the
products $z_\alpha z_\beta$ for $\alpha\leq \beta$, and we note that
there are $10$ such coefficients in each orthogonality condition.
This yields a total of $5k(k-1)$ constraints on the generators
$M^\alpha_{ab}$. }
\end{rem}

\subsection{ADHM construction of noncommutative
instantons}\label{section ADHM construction of noncommutative
instantons} We are ready for the construction of charge $k$
noncommutative bundles with instanton connections. As in previous
sections, we have the $(2k+2)\times k$ algebra-valued matrices
\begin{align*}
\sigma_z &= M^1 \otimes z_1+M^2 \otimes z_2+M^3 \otimes z_3+M^4
\otimes z_4, \\
\sigma_{J(z)} &=-M^1 \otimes z_2^* + M^2 \otimes z_1^* -M^3\otimes
z_4^* + M^4 \otimes z_3^*
\end{align*}
which, as already observed, have the properties $\sigma_{J(z)}^\star \sigma_z=0$
 and $\sigma_{J(z)}^\star\sigma_{J(z)}=\sigma_z^\star\sigma_z$.

\begin{lem}\label{lemma central rho generators} The entries of the
matrix $\rho^2:=\sigma_z^\star \sigma_z = \sigma_{J(z)}^\star\sigma_{J(z)}$
commute with those of
the matrix $\sigma_z$.
\end{lem}

\proof One finds that the $(\mu,\nu)$ entry of $\rho^2$ is
$$(\rho^2)_{\mu\nu}=\sum\nolimits_{r,\alpha,\beta}(M^\alpha{}^\dag)_{\mu
r}M^\beta_{r\nu}\otimes z_\alpha^*z_\beta$$ and that the $a,b$ entry of
$\sigma_z$ is $$(\sigma_z)_{ab}=\sum\nolimits_{\gamma} M^\gamma_{ab}\otimes
z_\gamma.$$ It is straightforward to check that these elements
always commute using the relations (\ref{eqn nc params}) for
$\A(\C^4_\theta)$ and the relations of Prop.~\ref{prop
universal matrix relations} for $\A(\tmM_\theta(H,K))$. The
essential feature is that every factor of $\eta_{\beta\alpha}$
coming from the relations between the $M^\alpha$'s is cancelled by a
factor of $\eta_{\alpha\beta}$ coming from the relations between the
$z_\alpha$'s.
\endproof

We need to enlarge slightly the matrix algebra $\M_k(\C)\otimes
\A(\tmM^{SD}_{\theta;k})\otimes \A(\C^4_\theta)$  by adjoining an
inverse element $\rho^{-2}$ for $\rho^2$, together with a square
root $\rho^{-1}$. That these matrices may be inverted is an
assumption, even in the commutative case where doing so corresponds
to the deletion of the non-generic points of the moduli space; these
correspond to so-called `instantons of zero size'.

From the previous lemma the matrix $\rho^2$, which is self-adjoint by construction, has entries in
the centre of the algebra $\A(\tmM^{SD}_{\theta;k})\otimes
\A(\C^4_\theta)$, so these new matrices $\rho^{-1}$ and $\rho^{-2}$
must also be self-adjoint with central entries. We collect the
matrices $\sigma_z$, $\sigma_{J(z)}$ together into the $(2k+2)\times
2k$ matrix
\begin{equation}\label{matU}
\V :=\begin{pmatrix}\sigma_z &
\sigma_{J(z)}\end{pmatrix}
\end{equation}
and we have by construction that
$$\V^* \V=\rho^2\begin{pmatrix} \mathbb{I}_{k} & 0 \\ 0 &
\mathbb{I}_k\end{pmatrix},$$ where $\mathbb{I}_k$ denotes the
$k\times k$ identity matrix. This of course means that the quantity
\begin{equation}\label{proj-q}
\Qp:=\V\rho^{-2} \V^*=\sigma_z\rho^{-2}\sigma_z^\star +
\sigma_{J(z)}\rho^{-2}\sigma_{J(z)}^\star
\end{equation}
is automatically a projection: $\Qp^2=\Qp=\Qp^*$ . For convenience we
denote
$$\Qp_z:=\sigma_z\rho^{-2}\sigma_z^\star, \qquad \Qp_{J(z)}:=\sigma_{J(z)}\rho^{-2}\sigma_{J(z)}^\star,$$
which are themselves projections, in fact orthogonal ones, $\Qp_{J(z)}
\Qp_z = 0$, due to the fact that $\sigma_{J(z)}^\star \sigma_z =0$.

\begin{lem} The trace of the projection $\Qp_z$ is equal to $k$; likewise for  $\Qp_{J(z)}$.\end{lem}

\proof We compute the trace as follows:
\begin{align*}\textup{Tr}                       %Achange at CMI: change in line breaks
\,\Qp_z&=\sum\nolimits_\mu (\sigma_z\rho^{-2}
\sigma_z^\star)_{\mu\mu}=\sum\nolimits_{\mu,r,s}(\sigma_z)_{\mu
r}(\rho^{-2})_{rs}(\sigma_z^\star)_{s\mu}\\
&=\sum\nolimits_{\mu,r,s}(\rho^{-2})_{rs}(\sigma_z)_{\mu
r}(\sigma_z)^\star_{s\mu}=\sum\nolimits_{\mu,r,s}(\rho^{-2})_{rs}(\sigma_z)^\star_{s\mu}(\sigma_z)_{\mu
r}\\
&=\sum\nolimits_{r,s}(\rho^{-2})_{rs}(\sigma_z^\star\sigma_z)_{sr}=
\textup{Tr}\, \mathbb{I}_k=k.
\end{align*} In the third equality we
have used the fact that, as said, the entries of $\rho^{-2}$ commute
with those of $\sigma_z$, whereas in the fourth equality we have
used the fact that every element of $\A(\tmM^{SD}_{\theta;k})\otimes
\A(\C^4_\theta)$ commutes with its own adjoint. An analogous chain
of equality establishes the same result for the projection
$\Qp_{J(z)}$.\endproof

As a consequence the projection $\Qp$ has trace $2k$.

\begin{prop}\label{prop ADHM projector} The operator
$$
\Pp:=\mathbb{I}_{2k+2}-\Qp
$$
is a projection in the algebra
$\M_{2k+2}\left(\A(\tmM^{SD}_{\theta;k})\otimes
\A(S^4_\theta)\right)$ with trace equal to $2$.\end{prop}

\proof The entries of the projection $\Qp_z$ are in the subalgebra
of $\A(\tmM^{SD}_{\theta;k})\otimes\A(\C^4_\theta)$ made up of
$\U(1)$-invariants which, by the discussion of Sect.~\ref{section
quantum conformal transformations}, is precisely
$\A(\tmM^{SD}_{\theta;k})\otimes \A(\CP^3_\theta)$. Now recall from
Sect.~\ref{section noncommutative twistor space} that the degree one
elements of $\A(\CP^3_\theta)$ of the form $Z+J(Z)$ generate the
$J$-invariant subalgebra, which may be identified with
$\A(S^4_\theta)$. The entries of $\Qp_z$ being linear in the
generators of $\A(\CP^3_\theta)$, it follows that the projection
$\Qp$ has entries in
$\A(\tmM^{SD}_{\theta;k})\otimes\A(S^4_\theta)$. The same is true of
the complementary projection $\Pp$ as well. Finally, since the
projection $\Qp$ has trace $2k$, the trace of the projector $\Pp$ is
just $2$.\endproof

We think of the projective right $\A(S^4_\theta)$-module
$\E:=\Pp\A(S^4_\theta)^{2k+2}$ as defining a family of rank two vector
bundles over $S^4_\theta$ parameterised by the noncommutative space
$\tmM^{SD}_{\theta;k}$. We equip this family of vector bundles with
the associated family of Grassmann connections $\n:=\Pp \circ
(\textup{id}\otimes \D)$, after extending the exterior derivative
from $\A(S^4_\theta)$ to $\A(\tmM^{SD}_{\theta;k})\otimes
\A(S^4_\theta)$ by $\textup{id}\otimes \D$. Moreover, we need also
to extend the Hodge $*$-operator as $\textup{id}\otimes *_\theta$.

\begin{prop} The curvature $F=\Pp((\textup{id}\otimes \D) \Pp)^2$ of the Grassmann connection $\n=\Pp\circ (\textup{id}\otimes \D)$ is anti-self-dual, that is to say
$(\textup{id}\otimes *_\theta) F = - F$.
\end{prop}

\proof When $\theta=0$ this construction is the usual ADHM
construction and it is known \cite{adhm:ci} ({\em cf}. also
\cite{mw:isdtt}) that it produces connections whose curvature is
an anti-self-dual two-form:
$$
(\textup{id}\otimes *_\theta)\, \Pp((\textup{id}\otimes \D) \Pp)^2=-\Pp((\textup{id}\otimes \D) \Pp)^2.
$$
As observed in Sect.~\ref{section noncommutative hopf fibration},
the Hodge $*$-operator is defined by the same formula as it is
classically and, as vector spaces, the self-dual and anti-self-dual
two-forms $\Omega^2_{\pm}(S^4_\theta)$ are the same as their
undeformed counterparts $\Omega^2_{\pm}(S^4)$. This identification
survives the tensoring by $\A(\tmM^{SD}_{\theta;k})$ which yields
$\A(\tmM^{SD}_{\theta;k})\otimes\Omega^2_{\pm}(S^4_\theta)$ to be
isomorphic, as vector spaces, to
$\A(\tmM^{SD}_{k})\otimes\Omega^2_{\pm}(S^4)$. Thus the
anti-self-duality holds also when $\theta \not=0$.
\endproof
\begin{rem}\textup{One may alternatively verify the anti-self-duality via a
complex structure. Indeed, there is an (almost) complex
structure $\gamma:\Omega^1(\CP^3_\theta) \rightarrow
\Omega^1(\CP^3_\theta)$  given by $\gamma(\D z_l):=(\D \circ J)(z_l)$,
$l=1,\ldots,4$, the operator $J$ being the one defined in \eqref{eqn J},
for which we declare the forms $\D z_l$ to be
holomorphic and the forms $\D z_l^*$ to be anti-holomorphic.
For instance, on
generators of $\A(\CP^3_\theta)$ we have
$$\D(z_jz_l^*)=\eta_{lj} z_l^*\D z_j + z_j\D z_l^*$$ from the Leibniz rule
and the relations \eqref{eqn nc params}, and we write
$\D=\partial + \bar \partial$ with respect to this decomposition
into holomorphic and anti-holomorphic forms. Since as vector spaces
the various graded components of the differential algebra
$\Omega(\CP^3_\theta)$ are undeformed, these operators $\partial$,
$\bar \partial$ extend to a full Dolbeault complex with
$\partial^2=\bar
\partial^2=\bar\partial\partial+\partial\bar\partial=0$.
The algebra inclusion $\A(S^4_\theta)\hookrightarrow
\A(\CP^3_\theta)$ extends to an inclusion of differential graded
algebras $\Omega(S^4_\theta)\hookrightarrow \Omega(\CP^3_\theta)$
and the Hodge operator $*_\theta$ is, as in the classical case,
defined in such a way that a two-form $\omega \in
\Omega^2(S^4_\theta)$ is anti-self-dual if and only if its image in
$\Omega^2(\CP^3_\theta)$ is of type $(1,1)$. Thus, to check that the
curvature $\Pp((\textup{id}\otimes \D)\Pp)^2$ is anti-self-dual, we
use this inclusion of forms (i.e. we express everything in terms of
$\D z_j$, $\D z_j^*$) and check that each component
$F_{ad}=\Pp_{ab}((\textup{id}\otimes \D)\Pp_{bc})\wedge
(\textup{id}\otimes \D)\Pp_{cd})$ of the curvature is a sum of terms
of type $(1,1)$. This approach to noncommutative twistor theory,
including a more explicit description of the noncommutative
Penrose-Ward Transform, will be discussed in more detail elsewhere.
%\cite{bm:ip}
}
\end{rem}

We next turn to the computation of the topological charge of the
family of bundles $\E:=\Pp\A(S^4_\theta)^{2k+2}$ given above. To this
end we observe that the matrix $\sigma_z$ has $k$ linearly
independent columns (since if not, it would not be injective) and
that the columns of $\sigma_{J(z)}$ are obtained from those of
$\sigma_z$ by applying the map $J$. Clearly we are free to rearrange
the columns of the matrix $\V$ (since this will not alter the class
of the projection $\Pp$), whence we may as well arrange them as
$$\V=\begin{pmatrix}
\sigma_z{}^{(1)}
&J(\sigma_z{}^{(1)})&\sigma_z{}^{(2)}&J(\sigma_z{}^{(2)})
&\cdots&\sigma_z{}^{(k)}&J(\sigma_z{}^{(k)})\end{pmatrix},$$ where
$\sigma_z{}^{(l)}$ denotes the $l$-th column of $\sigma_z$ and
$J(\sigma_z{}^{(l)})$ denotes the $l$-th column of $\sigma_{J(z)}$.
For fixed $l$, we denote the
entries of the column $\sigma_z{}^{(l)}$ (together with their
conjugates) by
$$w_\mu{}^{(l)}:=\sum\nolimits_{\alpha} M^\alpha_{\mu l}\otimes z_\alpha,
\quad (w_\mu{}^{(l)}){}^*:=\sum\nolimits_{\alpha} M^\alpha_{\mu
l}{}^*\otimes z_\alpha^*, \qquad \mu=1,\ldots,2k+2.$$ The entries of
the column $J(\sigma_z{}^{(l)})$ are obtained from those of
$\sigma_z{}^{(l)}$ by applying the map $J$, and one clearly has
$J((w_\mu{}^{(l)}){}^*)=(J(w_\mu{}^{(l)}))^*$. In the classical
limit $\theta=0$, we could evaluate the parameters $M^\alpha_{ab}$
as fixed numerical values: this would identify the columns
$\sigma_z^{(l)}$ and $J(\sigma_z^{(l)})$ as spanning a quaternionic
line in $\HH^{k+1}$, where the latter is defined by the $2k+2$
complex coordinates $w_\mu^{(l)}$ and their conjugates, equipped
with an anti-involution $J$. In the noncommutative case, although we
lack the evaluation of the parameters $M^\alpha_{ab}$, we continue
to interpret the columns $\sigma_z^{(l)}$ and $J(\sigma_z^{(l)})$ as
spanning a `one-dimensional' quaternionic line.

As already observed in Rem.~\ref{rem orth columns}, the columns of
$\sigma_z$ are orthogonal, as are the columns of $\sigma_{J(z)}$;
whence the rank $2k$ projection $\Qp$ in \eqref{proj-q} decomposes as a sum of
projections
$$
\Qp=\Qp_1 + \cdots + \Qp_k,
$$
where $\Qp_l:=\widetilde\Psi_i\widetilde\Psi_l{}^*$ is the rank two
projection defined by the $(2k+2) \times 2$ matrix
$\widetilde\Psi_l$ comprised of the columns $\sigma_z{}^{(l)}$ and
$J(\sigma_z{}^{(l)})$, appropriately normalised by $\rho^{-1}$.
Explicitly, this matrix is
$$
\widetilde \Psi_l:=\begin{pmatrix}\sum_{r,\alpha}(M^\alpha_{\mu
r}\otimes z_ \alpha)(\rho^{-1})_{rl}& \sum_{s,\beta}(M^\beta_{\mu
s}\otimes
J(z_\beta))(\rho^{-1})_{sl}\end{pmatrix}_{\mu=1,\ldots,2k+2},
$$
and a direct check yields $\widetilde \Psi_l{}^* \widetilde
\Psi_l=\mathbb{I}_{2}$ so that $\Qp_l$ is indeed a projection for each
$l=1,\ldots,k$. Hence the matrix $\V$ in \eqref{matU} has $2k$ columns which we
interpret as spanning $k$ quaternionic lines, with the same being
true of the normalised matrix $\V\rho^{-1}$. The computation of the
topological charge of the projection $\Qp$ therefore boils down to the
computation of the charge of each of the projections $\Qp_l$, for
$l=1,\ldots,k$.

\begin{lem} For each $l=1,\ldots,k$ the projection $\Qp_l$ is
Murray-von Neumann equivalent to the projection $1\otimes \qp$ in the
algebra $\M_{2k+2}(\A(\tmM^{SD}_{\theta;k})\otimes \A(S^4_\theta))$,
where $\qp$ is the basic projection defined in equation (\ref{eqn modified
basic instanton projector}). \end{lem}

\proof From equations \eqref{eqn modified partial isometry} and
(\ref{eqn modified basic instanton projector}) we know that
$\qp=\Psi\Psi^*$. Then, for each $l=1,\ldots,k$ define a partial
isometry $V_l$ in $\M_{2k+2,4}\left(\A(\tmM^{SD}_{\theta;k})\otimes
\A(S^4_\theta)\right)$ by
$$
V_l:=\widetilde \Psi_l \, (1\otimes\Psi^*), \qquad
V_l^*:=(1\otimes\Psi) \, \widetilde \Psi_l^*.
$$
Straightforward computations show  that $V_lV_l{}^*=\Qp_l$ and
$V_l{}^* V_l=1\otimes \qp$.
\endproof

We invoke the strategy of \cite{lprs:ncfi} to compute the
topological charge of the family of bundles defined by each $\Qp_l$.
Indeed, the charge of the projection $\qp$ was shown in
\cite{lvs:pfns} be equal to $1$, given as a pairing between the
second Chern class $\textup{ch}_2(\qp)$, which lives in the cyclic
homology group $\textup{HC}_4(\A(S^4_\theta))$, with the fundamental
class of $S^4_\theta$, which lives in the cyclic cohomology
$\textup{HC}^4(\A(S^4_\theta))$. Although the class
$\textup{ch}_2(\Qp_l)$, being an element in
$\textup{HC}_4(\A(\tmM^{SD}_{\theta;k})\otimes \A(S^4_\theta))$, may
not {\em a priori} be paired with the fundamental class of
$S^4_\theta$, Kasparov's KK-theory is used to show that in fact
there is a well-defined pairing between the K-theory
$\textup{K}_0\left(\A(\tmM^{SD}_{\theta;k})\otimes
\A(S^4_\theta)\right)$ and the K-homology
$\textup{K}^0(\A(S^4_\theta))$. Since by the previous lemma the
projections $1\otimes \qp$ and $\Qp_l$ define the same class in the
K-theory of $\A(\tmM^{SD}_{\theta;k})\otimes \A(S^4_\theta)$, it
follows as in \cite{lprs:ncfi} that the topological charge of each
projection $\Qp_l$ is equal to $1$.

\begin{prop} The family of bundles $\E=\Pp\A(S^4_\theta)^{2k+2}$ has topological charge equal to $-k$. \end{prop}

\proof By the argument given above, the projections $\Qp_l$ have
topological charge equal to $1$ for each $l=1,\ldots,k$. The
projection $\Qp=\Qp_1+\cdots+\Qp_k$ therefore has charge $k$, whence
$\Pp$ must have charge $-k$.\endproof

We finish this section by remarking that the construction given
above in the section has an interpretation in terms of `universal
connections', as described in \cite{ma:gymf}. As already said, the
classical quaternion vector space $\HH^{k+1}$ may be identified with
the complex vector space $\C^{2k+2}$ equipped with the quaternionic
structure $J$. Points of the Grassmannian manifold
$\textup{Gr}_k(\HH^{k+1})$ of quaternionic $k$-dimensional subspaces
of $\HH^{k+1}$ may thus be identified with $2k$-dimensional
subspaces of $\C^{2k+2}$ which are invariant under the involution
$J$. Following the general strategy of \cite{bm:qtt} for the
coordinatisation of Grassmannians, the algebra of functions on
$\textup{Gr}_k(\HH^{k+1})$ is given by functions taking values in
the set of rank $2k$ projectors $P=(P^\mu{}_\nu)$ on $\C^{2k+2}$
which are $J$-invariant, viz.
\begin{multline*}                           %Achange at CMI: formula break
\A(\textup{Gr}_k(\HH^{k+1})):= \\
\C\left[P^\mu{}_\nu ~|~ \sum\nolimits_\lambda P^\mu{}_\lambda P^\lambda{}_\nu=P^\mu{}_\nu,~
(P^\mu{}_\nu)^*=P^\nu{}_\mu,~ \sum\nolimits_\mu P^\mu{}_\mu =2k, ~ ~J(P^\mu{}_\nu)=P^\mu{}_\nu\right],
\end{multline*}
where $\mu,\nu=1,\ldots,2k+2$. In the classical case, when
$\theta=0$, the projection $\Qp$ in \eqref{proj-q} realises
$\A(S^4_{\theta=0})$ as a subalgebra of
$\A(\textup{Gr}_k(\HH^{k+1}))$, whence this construction should be
viewed as the dual of an embedding $S^4 \hookrightarrow
\textup{Gr}_k(\HH^{k+1})$, as given in \cite{ma:gymf}. We expect
that, in the deformed case, the projection $\Qp$ views
$\A(S^4_\theta)$ as a subalgebra of a suitably-deformed version of
$\A(\textup{Gr}_k(\HH^{k+1}))$. For fixed $k$, the set of monads is
bound to parameterise the set of such `algebra embeddings'.

\subsection{ADHM construction of charge one instantons} \label{section adhm charge one}
As a way of illustration we briefly verify that the above ADHM
construction of noncommutative families of instantons gives back the
family constructed in \cite{lprs:ncfi} when performed for the charge
one case.

The starting point is the basic instanton on $S^4_\theta$ described
in Sect.~\ref{section basic instanton} and which arises via a monad
construction as follows. The monad we consider is the sequence
\begin{equation} \label{eqn charge one monad}
\A(\C^4_\theta)(-1) \xrightarrow{\sigma_z} \C^4 \otimes
\A(\C^4_\theta) \xrightarrow{\tau_z} \A(\C^4_\theta)(1),
\end{equation}
where the arrows are the maps
$$\sigma_z=\begin{pmatrix}z_1&z_2&z_3&z_4\end{pmatrix}^{\textup{t}}, \qquad \tau_z=
\sigma_{J(z)}^\star=\begin{pmatrix}-z_2&z_1&-z_4&z_3\end{pmatrix}.$$
Since $\tau_z\sigma_z=\sigma_{J(z)}^\star\sigma_z=0$, it is clear
that this is a monad with $k=1$; by construction it is self-dual. In
the present case $\rho^2=\sigma_z^\star \sigma_z=\sum_j
z_j^*z_j=r^2$, which we already assumed was invertible
(corresponding to the deletion of the origin in $\C^4_\theta$). One
computes that
$$\V\V^*=\tfrac{1}{2}r^{-2}\begin{pmatrix}r^2+x&0&\alpha&\beta\\0&r^2+x&-\mu \beta^* &\bar \mu\alpha^*\\ \alpha^*
& -\bar\mu\beta & r^2-x&0\\ \beta^*&\mu\alpha&0&r^2-x\end{pmatrix}$$
which is just the projector $\qp$ of equation (\ref{eqn modified basic
instanton projector}). This is the `tautological' monad construction
given in \cite{bm:qtt}. The anti-self-dual version is the
projector $\Pp=1-\V\V^*$, in agreement with the ADHM construction
above.

The monad (\ref{eqn charge one monad}) may be rewritten in the form
\begin{equation}\label{eqn charge one sigma}
\sigma_z= (1, 0, 0, 0)^{\textup{t}} \otimes z_1~+~ (0,1,0,0)^{\textup{t}} \otimes z_2
~+~(0,0,1,0)^{\textup{t}} \otimes z_3~+~(0,0,0,1)^{\textup{t}} \otimes z_4,
\end{equation}
with $\tau_z$ defined as its dual. With the
strategy of \cite{lprs:ncfi} one generates new instantons by
coacting on the generators $z_1,\ldots,z_4$ with the quantum
conformal group $\A(\SL_\theta(2,\HH))$. Using the formula (\ref{eqn
conformal coaction}) for the coaction, the monad map (\ref{eqn
charge one sigma}) transforms into
\begin{equation}\sigma_{\Delta_L(z)}=\begin{pmatrix}a_1\\a_2\\c_1\\c_2\end{pmatrix}\otimes
z_1~+~\begin{pmatrix}-a_2^*\\a_1^*\\-c_2^*\\c_1^*\end{pmatrix}\otimes
z_2~+~\begin{pmatrix}b_1\\b_2\\d_1\\d_2\end{pmatrix}\otimes
z_3~+~\begin{pmatrix}-b_2^*\\b_1^*\\-d_2^*\\d_1^*\end{pmatrix}\otimes
z_4,\end{equation} and these four column vectors are the columns
of the matrix (\ref{eqn defining M(H)}) which defines the algebra
$\A(\SL_\theta(2,\HH))$. If we write
$$
\widehat{M}^{1}=\begin{pmatrix}a_1&a_2&c_1&c_2\end{pmatrix}^{\textup{t}},
\qquad
\widehat{M}^2=\begin{pmatrix}-a_2^*&a_1^*&-c_2^*&c_1^*\end{pmatrix}^{\textup{t}},
$$
$$
\widehat{M}^3=\begin{pmatrix}b_1&b_2&d_1&d_2\end{pmatrix}^{\textup{t}}, \qquad
\widehat{M}^4=\begin{pmatrix}-b_2^*&b_1^*&-d_2^*&d_1^*\end{pmatrix}^{\textup{t}},
$$
then we have the algebra relations $\widehat{M}^\alpha_j \widehat{M}^\beta_l =
\eta_{jl}\eta_{\beta \alpha}\widehat{M}^\beta_l \widehat{M}^\alpha_j$ coming
from the relations (\ref{eqn twisted quantum group relations}) for $\A(\SL_\theta(2,\HH))$. We
thus think of the algebra generated by the $\widehat{M}_j^\alpha$ as
parameterising the set of charge one instantons, since the map
$\sigma_{\Delta_L(z)}$ may be used to construct the family (\ref{eqn
coacted instanton projector}) of projections with topological charge
equal to $1$ and hence a family of Grassmann connections with
anti-self-dual curvature, just as in \cite{lprs:ncfi}.

In contrast, the ADHM construction of Sect.~\ref{section ADHM construction of noncommutative
instantons} for the case $k=1$ says that the
charge one monads are parameterised by the algebra
$\A(\tmM_{\theta;k})$ generated by the matrix elements $M_j^\alpha$,
with $j,\alpha=1,\ldots,4$, subject in particular to the relations
$M^\alpha_jM^\beta_l=\eta_{\beta \alpha}M^\beta_lM^\alpha_j$.

We see that these two approaches seem to give different
parameterisations of the set of monads for the case $k=1$, and hence
of the set of charge one instantons. The discrepancy has its root in
the fact that the ADHM construction requires generators lying in the
same row of the matrix $(A_{ij})$ to commute, whereas the `coaction
approach' given above says that such generators do not commute.

However, the discrepancy fades away when we pass to the `true' parameter space
for the families. On the one hand, as observed in \cite{lprs:ncfi},
the coaction \eqref{eqn unitary coaction} of the
quantum subgroup $\A(\Sp_\theta(2))$ of $\A(\SL_\theta(2,\HH))$ leaves the basic one-form
\eqref{eqn basic instanton one form} invariant.
We think of the
latter coaction as generating gauge-equivalent instantons, so that
the `true' parameter space for this family is rather the subalgebra
of $\A(\SL_\theta(2,\HH))$ of coinvariants under the coaction of
$\A(\Sp_\theta(2))$. The generators of this algebra are computed to
be
$$
\widehat{m}_{\alpha\beta}:=\sum\nolimits_l \widehat{M}_l^\alpha{}^* \widehat{M}_l^\beta,\qquad \alpha,\beta=1\ldots,4,
$$
whose relations are easily found to be
\begin{equation}\label{eqn moduli generators}
\widehat{m}_{\alpha\beta} \widehat{m}_{\mu\nu}=\eta_{\beta\mu}\eta_{\nu\beta}\eta_{\mu\alpha}\eta_{\alpha\nu}\widehat{m}_{\mu\nu}\widehat{m}_{\alpha\beta},
\end{equation}
and which certainly do not
depend on the rows of the matrix $(A_{ij})$. On the other hand,
gauge equivalence for the ADHM family parameterised by the
$M_i^\alpha$ is generated by the action of the classical group
$\Sp(2)$ (we borrow this result from Prop.~\ref{prop symplectic
gauge} in the next section), and here the invariant subalgebra is
generated by elements of the form
$$m_{\alpha\beta}:=\sum\nolimits_l M_l^\alpha{}^*M_l^\beta, \qquad
\alpha,\beta=1\ldots,4.$$ The relations in this algebra are just as
in equation \eqref{eqn moduli generators}, so that these two
families of charge one instantons are just the same.

\section{Gauge Equivalence of Noncommutative Instantons}
\label{section gauge theory}

Classically, a way to think of a gauge transformation of a vector
bundle $E$ over $S^4$ is as a unitary change of basis in each fibre
$E_x$ in a way which depends smoothly on $x \in S^4$.
Two connections on $E$ are said to be gauge equivalent if they are related by a
gauge transformation in this way. Now, rather than being interested
in the set of all instantons on $S^4$, one is interested in the
collection of gauge equivalence classes, that is to say classes of
instantons modulo gauge transformations.

It is therefore necessary to have an analogue of the notion of gauge
equivalence also for the noncommutative families of instantons
constructed previously. In fact, noncommutative geometry is a very
natural setting for the study of gauge transformations, as we shall
see in this section; we refer in particular to \cite{ac:fncg,ac:rec}
({\em cf}. also \cite{gl:book}).

\subsection{Gauge equivalence for families of instantons}\label{section usual gauge}

Recall that a first order differential calculus on a unital $*$-algebra $A$ is a pair
$(\Omega^1A,\D_A)$, where $\Omega^1A$ is an $A$-$A$-bimodule giving
the space of one-forms and $\D_A:A \to \Omega^1A$ is a linear map
satisfying the Leibniz rule,
$$\D_A (x y)=x(\D_A y)+(\D_A x)y \qquad \text{for all} ~~\quad x,y \in A.$$
One also assumes that the map $x \otimes y \to x(\D_A y)$ is
surjective. One names $\Omega^1A$ a $*$-calculus if for
$x_j,y_j\in A$ one has that $\sum\nolimits_j x_j\D y_j=0$ implies
$\sum\nolimits_j\D(y_j^*)x_j^*=0$: it follows from this condition
that there is \cite{wor} a unique $*$-structure on $\Omega^1A$ such
that $(\D_A a)^*=\D_A (a^*)$ for all $a \in A$. The differential
calculi on $\A(S^4_\theta)$ and $\A(S^7_\theta)$ in
Sect.~\ref{section noncommutative hopf fibration} are examples of
first order differential $*$-calculi on noncommutative spaces.

Let us fix a choice of $*$-calculus on $A$. Then let $\E$ be a
finitely generated projective right $A$-module endowed with an
$A$-valued Hermitian structure denoted by $\la \cdot |\cdot \ra$. A
connection on $\E$ is a linear map $\n:\E\to \E\otimes_A \Omega^1A$
satisfying the Leibniz rule
$$
\n(\xi x)=(\n \xi)x + \xi \otimes \D_A x \qquad \text{for all} ~~\xi
\in \E, ~x \in A.
$$
The connection $\n$ is said to be compatible with the Hermitian
structure on $\E$ if it obeys
$$\la \n \xi|\eta\ra + \la\xi|\n \eta
\ra= \D_A \la \xi|\eta \ra \qquad \text{for all} ~~\xi
\in \E, ~x \in A.
$$
On $\E$ there is at least one compatible connection, the so-called
Grassmann connection $\n_0$. If $P\in \End_A(\E)$, $P^2=P=P^*$, is
the projection which defines $\E$ as a direct summand of a free
module, that is, $\E=P(\C^N\otimes A)$, then $\n_0=P\circ \D$. Any
other connection on $\E$ is of the form $\n=\n_0+\omega$, where
$\omega$ is an element of $\textup{Hom}_A(\E,\E\otimes_A
\Omega^1A)$.

The gauge group of $\E$ is defined to be
$$
\mathcal{U}(\E):=\left\{ U \in \textup{End}_A(\E)~|~ \la
U\xi|U\eta\ra=\la\xi|\eta\ra ~\text{for all}~ \xi,\eta \in \E
\right\}.
$$
If $\n$ is a compatible connection on $\E$, each element $U$ of the
gauge group $\mathcal{U}(\E)$ induces a `new' connection by the action
$$\n^{U}:=U\n U^*.
$$ Of course, $\n^U$ is not really a different
connection, it simply expresses $\n$ in terms of the transformed
bundle $U\E$, hence one says that a pair of connections $\n_1$, $\n_2$
on $\E$ are {\em gauge equivalent} if they are related by such a
gauge transformation $U$. In terms of the decomposition
$\n=\n_0+\omega$, one finds that $\n^U= \n_0+\omega^U$, where
$$\omega^U:=U (\n_0 U^*) + U\omega U^*.$$
%For the Grassmann connection,
A choice of gauge would be a choice of partial isometry
$\Psi:\E\rightarrow \A^N$ such that $\Psi^*\Psi=\textrm{Id}_{\E}$
and $\Psi\Psi^*=P$. Any other gauge is then given by an element $U$
of the gauge group of $\E$: the partial isometry $\Psi$ gets
replaced by $U\Psi$, for which we indeed have
$$(U\Psi)^*(U\Psi)=\Psi^*\Psi=\textrm{Id}_{\E}.$$
and the projection $P$ gets transformed to
$$(U\Psi)(U\Psi)^*=U(\Psi\Psi^*)U^*=UPU^*,$$
an operation that does not change the equivalence class of $P$. In
the fixed gauge the Grassmann connection $\n_0=P\circ \D$ naturally acts on
`equivariant maps' $\varphi=\Psi F$ where $F\in \A^N$. The result is an `equivariant one-form',
$$
\n_0(\Psi F)= (\Psi\Psi^*) \D (\Psi F) = \Psi \Big( \D F +  \Psi^* \D( \Psi) F \Big),
$$
and identifies the \emph{gauge potential} to be given by
$$
A= \tfrac{1}{2} \left( \Psi^* (\D\Psi) - \right  (\D \Psi^*) \Psi ).
$$
Under the transformation $\Psi \mapsto U\Psi$, the gauge potential transforms as expected:
 $$
\Psi^* \D \Psi \mapsto \Psi^*(\D\Psi) + \Psi^*U^*(\D U)\Psi.
$$

We now turn back to the construction of instantons. Gauge
equivalence being defined as above by unitary module endomorphisms
means that we are free to act on the right $\A(\C^4_\theta)$-module
$\mathcal{K}=K\otimes \A(\C^4_\theta)$ by a unitary element of the
matrix algebra $\M_{2k+2}(\C)\otimes \A(\C^4_\theta)$. In order to
preserve the instanton construction, we must do so in a way
preserving the bilinear form $(\,\cdot\,,\,\cdot\,)$ of equation
(\ref{eqn bilinear form}) which comes from the identification of $K$
with its dual $K^*$. Hence the map $\sigma_z$ in \eqref{eqn nc
module map} (or in \eqref{eqn dual nc module map}) is defined up to
a transformation $A \in
\textup{End}_{\A(\C^4_\theta)}(\mathcal{K})$, which is unitary and
is required to commute with the quaternion structure $J$. Similarly,
we are free to change basis in the modules $\mathcal{H}=H\otimes
\A(\C^4_\theta)$ and $\mathcal{L}=L\otimes\A(\C^4_\theta)$, provided
we preserve the fact that we identify $J(\mathcal{H})^\star$ and
$\mathcal{L}$. This means that the map $\tau_z$ of \eqref{eqn nc tau
expression} is defined up to an invertible transformation $B\in
\textup{End}_{\A(\C^4_\theta)}(\mathcal{H})$.

All this is saying is that the monad maps $\sigma_z$ and $\tau_z$ were
expressed as matrices with respect to a choice of basis for each of
the vector spaces $H$, $K$ and $L$; it is natural to question the
extent to which the resulting Grassmann connection $\n$ depends on
the choice of these bases. We denote by $\GL(\mathcal{H})$ the set
of automorphisms of $\mathcal{H}$ and by $\textup{Sp}(\mathcal{K})$
the set of all unitary endomorphisms of $\mathcal{K}$ respecting the
quaternion structure:
$$\textup{Sp}(\mathcal{K}):=\{ A \in
\textup{End}_{\A(\C^4_\theta)}(\mathcal{K})~|~\la A\xi|A\xi\ra=\la
\xi|\xi\ra,~J(A\xi)=AJ(\xi)~\text{for all}~\xi \in \mathcal{K}\}.$$
Given $A \in \textup{Sp}(\mathcal{K})$ and $B \in \GL(\mathcal{H})$,
the gauge freedom is to map $\sigma_z \mapsto A\sigma_zB$.

\begin{prop} For all $B\in\GL(\mathcal{H})$, under the transformation
$\sigma_z\mapsto\sigma_zB$ the projection $\Pp$ of Prop.~\ref{prop ADHM projector} is left
invariant.
\end{prop}

\proof One first checks that $\rho^2\mapsto (\sigma_zB)^\star
(\sigma_zB) = B^\star \rho^2 B$ under this transformation, so that
$$\Qp_z \mapsto \sigma_z B (B^\star \rho^2 B)^{-1}B^\star
\sigma_z^\star = \sigma_z B (B^{-1} \rho^{-2} (B^\star)^{-1})
B^\star \sigma_z^\star = \Qp_z,$$ whence the projection $\Pp$ is
unchanged.\endproof

\begin{prop}\label{prop symplectic gauge}For all $A\in \textup{Sp}(\mathcal{K})$, under the transformation
$\sigma_z \mapsto A\sigma_z$ the projection $\Pp$ of Prop.~\ref{prop ADHM projector} transforms as
$\Pp\mapsto A\Pp A^\star$.
\end{prop}

\proof Replacing $\sigma_z$ by $A \sigma_z$ leaves $\rho^2$
invariant (since $A$ is unitary) and so has the effect that
$$\Qp_z \mapsto A\sigma_z \rho^{-2} \sigma_z^\star A^\star =
A\Qp_zA^\star,$$ whence it follows that $\Pp$ is mapped to
$A\Pp A^\star$.\endproof

These results give the general gauge freedom on monads, although
from the point of view of computing the number of constraints on the
algebra generators $M^\alpha_{ab}$  we need only consider the effect
of these transformations on the vector spaces $H$, $K$ and $L$, i.e.
it is enough to consider the groups of `constant' automorphisms.
This means the group $\Sp(K)=\Sp(k+1) \subset \Sp(\mathcal{K})$ and
the group $\GL(k,\RR)\subset \GL(\mathcal{H})$ (the latter because
we must preserve the identification of $J(\mathcal{H})^\star$ with
$\mathcal{L}$, and complex linear transformations of $H$ would
interfere with the tensor product in $J(\mathcal{H})^\star$). In
fact, it is known in the classical case that these constant
transformations are sufficient to generate all gauge symmetries of
the instanton bundles produced by the ADHM construction.

We conclude that in the noncommutative case as well
the gauge equivalence imposes an additional
$(k+1)(2(k+1)+1)$ constraints due to $\textup{Sp}(k+1)$ and a
further $k^2$ constraints due to $\textup{GL}(k,\RR)$. From
Rem.~\ref{rem number of params}, the total number of generators
minus the total number of constraints is thus computed to be
$$(8k^2+8k)-5k(k-1)-(3k^2+5k+3)=8k-3,$$ just as for the classical
case, a result which is somehow reassuring.

\subsection{Morita equivalent geometries and gauge theory}
\label{section gauge theory in general} It is a known idea that
Morita equivalent algebras describe the same topological space. The
simplest case is that of a one-point space $X=\{*\}$: the matrix
algebras $\M_n(\C)$ for any positive integer $n$ all have the same
one-point spectrum. More generally, if $X$ is a compact Hausdorff
space, the algebras $C(X)\otimes \M_n(\C)$ are all Morita equivalent
and all have the same spectrum $X$.

With this in mind, gauge theory arises naturally out of the
consideration of how to transfer differential structures between
Morita equivalent algebras. If one takes such structures to be
defined by a Dirac operator and associated spectral triple, then
the method for doing this is discussed in \cite{ac:fncg,ac:rec}.
Here we discuss a more general framework, where algebras may be
equipped with differential calculi not necessarily coming from a
spectral triple.

Let $A$ be a unital $*$-algebra and suppose that the $*$-algebra $B$
is Morita equivalent to $A$ {\em via} the $B$-$A$-bimodule $\E$,
that is to say $B \simeq \textup{End}_A(\E)$. In addition, on $\E$
there are compatible $A$-valued and $B$-valued Hermitian structures
\footnote{We shall also require the Hermitian structures to be
self-dual, {\em i.e.} every right $A$-module homomorphism $\varphi:
\E \to A$ is represented by an element of $\eta \in \E$ by the
assignment $\varphi(\cdot) = \la \eta | \cdot \ra$. A similar
property holds for the second Hermitian structure as well.}. Then a
choice of a connection $\nabla$ on $\E$, viewed as a right
$A$-module, yields a differential calculus on $B$. First of all, the
operator on $B$ given by
$$
\D_B^\n(x):=[\n,x], \qquad x \in B,
$$
is easily seen to be a derivation: $\D_B^\n (x y)=x(\D_B^\n
y)+(\D_B^\n x)y$, for $x,y \in B$. The $B$-$B$-bimodule $\Omega^1B$
of one-forms is then defined by
$$
\Omega^1B:= B \left( \D_B^\n(B) \right) B .
$$
For this to define a $*$-calculus we need that the connection $\n$
be compatible with the $A$-valued Hermitian structure on
$\E$ in the sense that
$$\la \n \xi|\eta\ra + \la\xi|\n \eta \ra=
\D_A \la \xi|\eta \ra
$$
for all $a \in A$ and $\xi,\eta \in \E$. If this compatibility condition is
satisfied, the assumption $\sum\nolimits_j x_j\D_B^\n y_j=0$
translates into $\sum\nolimits_j(x_j\n y_j)\xi=\sum\nolimits_j
x_jy_j(\n\xi)$ for all $x_j,y_j\in B$ and all $\xi \in \E$. This
implies, for all $\xi,\eta \in \E$ and all $x_j,y_j\in B$, that
\begin{align*}\sum\nolimits_j\la\D_B^\n(y_j^*)x_j^*\xi|\eta\ra&=\sum\nolimits_j\la\n(x_j^*y_j^*\xi)-x_j^*\n(y_j^*\xi)|\eta\ra \\
&=\sum\nolimits_j-\la x_j^*y_j^*\xi|\n\eta\ra+\D_A\la x_j^*y_j^*\xi|\eta\ra \\
& \qquad \qquad \qquad \qquad \qquad +\la y_j^*\xi|\n(x_j\eta)\ra-\D_A\la y_j^*\xi|x_j\eta\ra \\            %Achange at CMI: line break
&=\sum\nolimits_j-\la x_j^*y_j^*\xi|\n\eta\ra
+\la\xi|y_jx_j\n\eta\ra,\end{align*} whence it follows that
$\sum\nolimits_j\D_B^\n(y_j^*)x_j^*=0$ as it should for a
$*$-calculus. We interpret the passage $\D_A \to \D_B^\n$ as an {\em
inner fluctuation} of the geometry which results in a `Morita
equivalent' first order calculus $(\Omega^1B, \D_B^\n)$, now for the
algebra $B$.

A natural application is to think of the algebra $A$ as
being Morita equivalent to itself, so that $\E=A$ as a right
$A$-module and $B=A$. In this case, any Hermitian connection on $\E$
is necessarily of the form
\begin{equation}\label{herconn}
\n\xi=\D_A\xi + \omega\xi, \qquad \textup{for} \quad \xi \in\E,
\end{equation}
with $\omega=-\omega^* \in \Omega^1A$ a skew-adjoint one-form. The
corresponding differential on $B=A$ is computed to be
$$
(\D_A^\n  b) \xi = [\n,b]\xi=\n(b\xi)-b\n \xi=\D_A(b\xi)+\omega b\xi
-b\D_A\xi - b\omega\xi = (\D_A b)\xi + [\omega,b]\xi,
$$
using the Leibniz rule for $\D_A$. The passage
$$
\D_A \to \D_A^\n = \D_A + [\omega, ~\cdot~]
$$
is once again interpreted as an inner fluctuation of the geometry,
although when $A$ is commutative there are no non-trivial inner
fluctuations and thus no new degrees of freedom generated by the
above self-Morita mechanism. However, in the noncommutative
situation there is an interesting special case where $\omega$ is
taken to be of the form $\omega=u^*\D_A u$, for $u$ a unitary
element of the algebra $A$. Such a fluctuation is unitarily
equivalent to acting on $A$ by the inner automorphism
$$\alpha_u:A\rightarrow A, \qquad \alpha_u(a)=uau^*,$$ since for all $a \in A$ we
have that $\D_A^\n(a)=u^*\D_A(\alpha_u(a))u$. It therefore follows
that inner fluctuations defined by inner automorphisms generate
gauge theory on $A$.

\subsection{Gauge theory from quantum symmetries} We now consider a
slightly different type of gauge equivalence for our instanton
construction which is not present in the classical case and is a
purely quantum ({\em i.e.} noncommutative) phenomenon.

We consider the case where $A$ is a comodule $*$-algebra under a
left coaction of a Hopf algebra $H$, so that $A$ is isomorphic to
its image $B=\Delta_L(A)$. To transfer a calculus on $A$ to one on
$B$, a possible strategy is as follows. We take the $B$-$A$-bimodule
to be $\E:=B=\Delta_L(A)$ with left $B$-action and right
$A$-action defined by
$$b\tr\xi:=b\xi, \qquad \xi \tl a = \xi \Delta_L(a)$$ for $\xi \in
\E$, $a\in A$, $b \in B$. We also assume that the calculus
$\Omega^1A$ is left $H$-covariant, so that $\Delta_L$ extends to a
coaction on $\Omega^1A$ as a bimodule map such that $\D_A$ is an
intertwiner, whence the above bimodule structure on $\E$ extends to
one-forms in the natural way. This also canonically equips $B$ with
a $*$-calculus $\Omega^1B$, where the differential is
$\D_B=\textup{id}\otimes \D_A$.

We choose an arbitrary Hermitian connection on the right $A$-module
$\E$ for the calculus $(\Omega^1A,\D_A)$, which is necessarily of
the form
$$\n\xi=(\textup{id}\otimes\D_A)\xi + \tilde\omega\xi, \qquad \xi \in\E$$ with
$\tilde\omega=\Delta_L(\omega)$ for some $\omega=-\omega^* \in
\Omega^1A$ a skew-adjoint one-form. The corresponding differential
on $B$ is again defined by
$$
(\D_B^\n b) \xi = [\n, b\tr]\xi=\n(b\tr\xi)-b\tr\n
\xi=\D_A(b\tr\xi)+\omega (b\tr\xi) -b\tr \left( \D_A\xi+\omega\xi
\right),$$
and works out to be
$$\D_Bb=(\textup{id}\otimes \D_A)b + [\tilde \omega,b].$$ Note also that for all $b\in B$ we have $b=\Delta_L(a)$ for some $a
\in A$ and so it follows that
$$\D_Bb=\Delta_L(\D_Aa) + \Delta_L([\omega,a]),$$ so that the
coaction commutes with inner fluctuations. Moreover, in the case
where $A$ is noncommutative, there are non-trivial inner
automorphisms of $A$ and hence
non-trivial gauge degrees of freedom
which
carry over from $A$ to $\Delta_L(A)$.

In particular, we apply this to the case $A=\A(S^4_\theta)$, with
$H=\A(\SL_\theta(2,\HH))$ the quantum conformal group of
$S^4_\theta $. The above discussion means that the coaction of
$\A(\SL_\theta(2,\HH))$ on $\A(S^4_\theta)$ by conformal
transformations in itself generates gauge freedom. The natural way
to extend the exterior derivative $\D_A$ on $\A(S^4_\theta)$ to
$\Delta_L(\A(S^4_\theta))$ is as $\textup{id}\otimes \D_A$: this
corresponds to taking $\tilde \omega=0$ and is the choice made in
\cite{lprs:ncfi}. However, in general we have the freedom to make
the transition $$\D_A \to (\textup{id}\otimes\D_A) + [\tilde \omega,
~\cdot~]$$ for some $\tilde \omega=\Delta_L(u^*\D_A u)$, where $u$
is some unitary element of $\A(S^4_\theta)$. Since the group of
inner automorphisms of $A$ is trivial when $A$ is commutative, this
is a feature of gauge theory which is certainly not present in the
classical case and is unique to the noncommutative paradigm.
More on this will be reported elsewhere.

\end{document}